\documentclass[12pt]{article}

\usepackage{amsmath,amssymb}
\usepackage{psfrag,graphicx}
\usepackage{verbatim}
\usepackage{color}

\usepackage[utf8]{inputenc}
\usepackage{verbatim}
\usepackage{subcaption}
\usepackage{drawmatrix}
\usepackage[thinlines]{easybmat}

\usepackage[a4paper]{geometry}
\newcommand{\fv}{\mathbf f}
\newcommand{\ev}{\mathbf e}
\newcommand{\xv}{\mathbf x}
\newcommand{\Xv}{\mathbf X}
\newcommand{\av}{\mathbf a}
\newcommand{\bv}{\mathbf b}
\newcommand{\uv}{\mathbf{u}}
\newcommand{\vv}{\mathbf v}
\newcommand{\nv}{\mathbf n}

\newcommand{\Ev}{ \mathbf{E}}
\newcommand{\nullv}{\mathbf{0}}

\newcommand{\tauv}{\boldsymbol{\tau}}
\newcommand{\Phiv}{\boldsymbol{\Phi}}
\newcommand{\Dv}{ \mathbf{D}}

\newcommand{\lamv}{{\lambda}}
\newcommand{\qv}{ \mathbf{q}}

\newcommand{\sigmat}{\boldsymbol{\sigma}}
\newcommand{\taut}{\boldsymbol{\tau}}
\newcommand{\eps}{\varepsilon}
\newcommand{\epst}{\boldsymbol{\varepsilon}}
\newcommand{\epsilont}{\boldsymbol{\epsilon}}
\newcommand{\St}{\mathbf{S}}
\newcommand{\Ct}{\mathbf{C}}

\newcommand{\dt}{\mathbf{d}}
\newcommand{\et}{\mathbf{e}}

\newcommand{\Ft}{\mathbf{F}}
\newcommand{\Ht}{\mathbf{H}}
\newcommand{\It}{\mathbf{I}}
\newcommand{\Et}{\mathbf{E}}

\newcommand{\Nt}{\mathbf{N}}
\newcommand{\at}{\mathbf{a}}
\newcommand{\bt}{\mathbf{b}}
\newcommand{\nullt}{\mathbf{0}}

\newcommand{\Am}{\mathbf{A}}
\newcommand{\Mm}{\mathbf{M}}

\newcommand{\Dm}{\mathbf{D}}
\newcommand{\Bm}{\mathbf{B}}
\newcommand{\Cm}{\mathbf{C}}
\newcommand{\Em}{\mathbf{E}}
\newcommand{\phim}{\boldsymbol{\varphi}}

\newcommand{\spaceV}{V}

\newcommand{\spaceSigma}{\Sigma}

\newcommand{\opdiv}{\operatorname{div}}

\graphicspath{
	{./}{./pics/}}
	
\title{3D mixed finite elements for curved, flat piezoelectric structures}
\author{Martin Meindlhumer\footnote{Johannes Kepler University Linz, Institute of Technical Mechanics, Altenbergerstr. 69, 4020 Linz, Austria,  e-mail: martin.meindlhumer@jku.at} and Astrid Pechstein\footnote{Johannes Kepler University Linz, Institute of Technical Mechanics, Altenbergerstr. 69, 4020 Linz, Austria, e-mail: astrid.pechstein@jku.at}}
	
\begin{document}

\maketitle

\begin{abstract}The Tangential-Displacement Normal-Normal-Stress (TDNNS) method is a finite element method that was originally introduced for elastic solids and later extended to piezoelectric materials. It uses tangential components of the displacement and normal components of the normal stress vector as degrees of freedom for elasticity. For the electric field, the electric potential is used. The TDNNS method has been shown to provide elements which do not suffer from shear locking. Therefore thin structures (e.g. piezoelectric patch actuators) can be modeled efficiently. 
	Hexahedral and prismatic elements of arbitrary polynomial order are provided in the current contribution. We show that these elements can be used to discretize curved, shell-like geometries by curved elements of high aspect ratio. The order of geometry approximation can be chosen independently from the polynomial order of the shape functions. We present two examples of curved geometries, a circular patch actor and a radially polarized piezoelectric semi-cylinder. Simulation results of the TDNNS method are compared to results gained in ABAQUS. 
	We obtain good results for displacements and electric potential as well as for stresses, strains and electric field when using only one element in thickness direction. 		
\end{abstract}

\paragraph{Keywords:} piezoelasticity, mixed finite elements, curved elements, locking free

\section{Introduction}

The direct and reverse piezoelectric effects allow to transform mechanical energy into electrical energy and vice versa. Moreover, modern piezoceramics are cheap in production and offer fast and accurate response. Therefore, piezoelectric materials are the method of choice in the realization of smart structures. 

When it comes to design and control of piezoelectric structures, simulation results are of high interest. A popular method to get approximate solutions is the finite element (FE) method. 
Allik and Huges \cite{AllikHughes:1970} and later Lerch \cite{Lerch:1988} carried out first simulations using volume elements based on the principle of virtual work. These elements are an extension of standard small-strain mechanical elements by degrees of freedom for the electric potential. However, piezoelectric structures are typically flat (e.g. patch actors), and therefore in order to avoid locking effects, their discretization into well-shaped elements leads to a computational system with a huge number of unknowns - even for simple applications. \par 

Nowadays, there are two ways to avoid this problem: the first is the derivation of equations for layered beams, plates and shells. The second is the design of locking-free volume elements. In any way, some relaxation technique is necessary to avoid shear locking. Prominent methods of relaxation are reduced integration, assumed strains, or mixed methods involving stresses, strains or dielectric displacements as additional independent unknowns. Also, a higher order of polynomial approximation is often feasible.

Mixed methods are probably the most involved of those methods mentioned above. Both from a theoretical point of view, when proving uniqueness and convergence results, as well as from an implementational point of view, when additional unknowns have to be realized in a computational code, the mixed methods require intricate treatment. However, when correctly done, mixed methods yield most accurate results, as contributions from various authors demonstrate. 
Among the mixed methods, we cite hybrid stress elements \cite{SzePan:1999,SzeYaoYi:2000}, a solid shell element based on a Hu-Washizu formulation \cite{KlinkelWagner:2006} or the multi variable formulation \cite{OrtigosaGil:2016}. Reissner-type mixed zigzag formulations were successfully used by \cite{CarreraBoscolo:2007,CarreraNali:2009,WuLin:2014}. High-order methods comprise, among many others, isogeometric and hierarchical elements from the group around Gabbert \cite{WillbergGabbert:2012,DuczekGabbert:2013}. A high-order element with zigzag function is proposed by Polit et al.~\cite{PolitDOttavioVidal:2016}.

  Pechstein and Sch{\"o}berl \cite{PechsteinSchoeberl:11} introduced an arbitrary order mixed FE-method for linear elasticity based on a Hellinger-Reissner formulation, which uses tangential displacements and  normal normal stresses as degrees of freedom (TDNNS-method). This method was shown to be locking-free \cite{PechsteinSchoeberl:12}, and later extended to piezoelectric materials \cite{PechsteinMeindlhumerHumer:2018}. In the current contribution, we want to show the capability of the method. We discuss high-order geometry approximation,
  such that curved geometries can be treated accurately by a low number of curvilinear elements. Moreover, we address the problem of eigenvalue computation for the mixed method. We propose to use an inverse iteration \cite{KNYAZEV2003} for the coupled electromechanical problem.

	%The order of geometry approximation is independent of the order of shape functions. Therefore curved geometries can be represented exactly. In our contribution we show the example of a circular patch actor. Furthermore we show how (possibly large) predeformations of piezoelectric structures can be considered. 

This paper is organized as follows: below, we introduce some notation on tensor products. In Section~\ref{sec:linear_piezoelasticity}, the problem of linear piezolelasticity is introduced. Section~\ref{sec:TDNNS_piezoelasticity} provides a short introduction into the variational formulation of the mixed TDNNS method that is the foundation of the current contribution. The implementation of the underlying shape functions and their derivatives is treated in Section~\ref{sec:curvilinear_elements}. The coupled eigenvalue problem is described in Section~\ref{sec:eigenvalues}, and an inverse iteration is proposed for its solution. Finally, numerical results illustrating the capability of the method are presented in Section~\ref{sec:num_results}.

\paragraph{Notation}
In the current contribution, we use tensor calculus. Vectorial and tensorial quantities are denoted by boldface symbols. When summing over tensor components we use Einstein's summation convention. The dot denotes contraction, while $\otimes^s$ is the symmetric tensor product:
\begin{align}
(\at \cdot \bt)_{ij} &= a_{ik} b_{kj}, & \at : \bt&= a_{ij} b_{ji}, & (\av \otimes^{s} \bv)_{ij} &= \tfrac{1}{2} ( a_i b_j + a_j b_i)
\end{align}
The gradient of a scalar field $a$ and a vector field $\bv$ is denoted by the nabla operator. We indicate with respect to which coordinates the differentiation is carried out,
\begin{align}
(\nabla_\xv a)_i &= \frac{\partial a}{\partial x_i}, & (\nabla_\xv \bv)_{ij} &= \frac{\partial b_i}{\partial x_j}.
\end{align} 
We may omit the index, if the derivation is with respect to $\xv$. The divergence operator of a vector field is defined in the standard way, for a tensor field it is understood row-wise,
\begin{align}
\opdiv_\xv \av &= \frac{\partial a_i}{\partial x_i}, & (\opdiv_\xv \bt)_i &= \frac{\partial b_{ij}}{\partial x_j}.
\end{align}

\section{Linear piezoelasticity} \label{sec:linear_piezoelasticity}

%\todo{Martin: check: alle Tensoren/Vektoren verwenden die Makros, alle nabla und div haben ein x, xhat, ... als index. Dann können wir später auch auf Vektorpfeile umstellen falls notwendig.}

We consider an elastic in part piezoelectric solid body in the domain $\Omega \subset \mathbb{R}^3$, subjected to body forces $\fv$ and free of body charges.We assume to stay in the regime of small deformations, where Voigt's theory of linear piezoelasticity \cite{voigt1910lehrbuch} is appropriate. We are interested in the displacement field  $\uv$ and in the electric potential $\phi$. From these quantities the electric field $\Ev=-\nabla\phi$ and the linear strain tensor $\epst=\frac{1}{2}(\nabla \uv+ {\nabla\uv}^T)$ are derived. For the static case the mechanical and electrical balance equations in differential form are given by
\begin{eqnarray}
\opdiv \sigmat &=& \fv, \label{eq:balance1}\\
\opdiv \Dv &=& 0. \label{eq:balance2}
\end{eqnarray}
with the Cauchy stress tensor $\sigmat$ and the dielectric displacement $\Dv$.  The mechanical boundary conditions are
\begin{equation}
\uv = \nullv \text{ on } \Gamma_1 \qquad \text{and} \qquad \sigmat \cdot \mathbf{n} = 0 \text { on } \Gamma_2 = \partial \Omega \backslash \Gamma_1,
\end{equation}
and the electrical boundary conditions
\begin{equation}
\phi = \phi_0 \text{ on } \Gamma_3 \qquad \text{and} \qquad \Dv \cdot \nv = \mathbf{0} \text { on } \Gamma_4 = \partial \Omega \backslash \Gamma_3.
\end{equation}

%\subsection{Constitutive Equations}

%\todo{Martin: alle $\backslash$ sind bei den Tensorprodukten eingefügt und alle transponierten stimmen..}

We are concerned with the most common variant of constitutive equations, which are often referred to as Voigt's linear theory of piezoelectricity \cite{voigt1910lehrbuch}. %\todo{Astrid oder Martin - Voigt zitat heraussuchen und einfügen}. 
In our notation, they read
\begin{eqnarray}
\sigmat &=& \Ct^E: \epst - \et^T\cdot \Ev, \label{eq:etype1}\\
\Dv &=& \et: \epst + \epsilont^{\eps}\cdot \Ev. \label{eq:etype2}
\end{eqnarray}
with $\Ct^E$ the elasticity tensor measured under constant electric field,  $ \epsilont^{\eps}$ the dielectric tensor at constant strain and the piezoelectric permittivity tensor $\et$. We use the following alternative variant of the material law %given in (\ref{eq:dtype1})- (\ref{eq:dtype2}).
\begin{align}
\epst &=\ \St^E: \sigmat + \dt^T\cdot \Ev, \label{eq:dtype1}\\
\Dv &=\ \dt: \sigmat + \epsilont^{\sigma}\cdot \Ev. \label{eq:dtype2}
\end{align}
Of course the material parameters are connected, their relations are %given by
\begin{align}
	 \St^{E} &=\ (\Ct^{E})^{-1},& \dt &=\ \et\cdot \St^{E}, & \epsilont^{\sigma} &=\ \epsilont^{\eps} + \dt\cdot \et^T.
\end{align}

\section{TDNNS for piezoelasticity} \label{sec:TDNNS_piezoelasticity}
%\todo{Astrid, section 3: Freiheitsgrad einführen, ausbauen}
In this section we briefly introduce the "Tangential Displacement Normal Normal Stress" (TDNNS) finite element method.  The TDNNS finite element method was introduced by Pechstein and Sch{\"o}berl \cite{PechsteinSchoeberl:11,pechstein2018analysis}   for  elasticity and later extended to piezoelasticity \cite{PechsteinMeindlhumerHumer:2018}. It is a mixed finite element method based on Reissner's principle which considers displacements and stresses as unknowns. A more detailed description can be found in \cite{PechsteinMeindlhumerHumer:2018}.
\par 
 We assume $\mathcal{T} = \{T\}$ to be a finite element mesh  of the domain $\Omega$. The unit vector $\nv$ denotes the outer normal on element or domain boundaries $\partial T$ or $\partial \Omega$. In general on a (boundary) surface a vector field $\vv$ can be split into a normal component $v_n=\vv\cdot\nv$ and a tangential component $\vv_t=\vv-v_n\nv$. A tensor field $\taut$ has a normal vector on a surface $\tauv_n=\taut\cdot\nv$, which analogously can be split into a normal component $\tau_{nn}$ and tangential component ${\tauv_{nt}}$.
 
\subsection{Finite element formulation for elastic problem } 
The TDNNS finite element method uses the tangential component of the displacement $\uv_t$ and the normal component of the stress vector $\sigma_{nn}$ as degrees of freedom. Those quantities are continuous across the element interfaces. Note that the normal displacement $u_n$ is discontinuous. For the displacements, tangential-continuous elements introduced by N{\'e}d{\'e}lec \cite{Nedelec:86} are used. Elements for the stresses are normal-normal-continuous and were introduced by Pechstein and Sch{\"o}berl \cite{PechsteinSchoeberl:11}. 
{For the electric potential, standard continuous (e.g. nodal or hierarchical) finite elements are used. }
We assume the tangential continuous Nedelec elements for the displacements and continuous nodal or hierarchical elements for the electric potential on hybrid meshes to be well known. The definition of arbitrary-order prismatic or hexahedral normal-normal continuous stress elements is postponed to Section \ref{sec:curvilinear_elements}. Other element types were introduced in \cite{PechsteinSchoeberl:11}. At present, it is only necessary to assume that admissible stresses $\sigmat$ are such that $\sigma_{nn}$ is continuous across all element interfaces.
\begin{comment}
The finite element space on simplicial meshes can be described by

\begin{eqnarray}\label{eq:fespace}
\uv, \delta \uv &\in& \spaceV_h = \{ \vv: \vv|_T \in [P^k(T)]^3, \vv_t \text{ continuous}\},  \label{eq:spaceVh}\\
\sigmat, \delta \sigmat &\in& \spaceSigma_h = \{ \taut: \taut|_T \in [P^k(T)]^{3\times3}_{sym}, \tau_{nn} \text{ continuous} \},\\  \label{eq:spaceSigmah}
\phi, \delta \phi &\in& \ \{ \phi:\phi_h |_T \in P^{k+1}(T),\phi \text{ continuous}\}\label{eq:wh}
%\text{with }  \spaceW_h &:=&\ \{ w: w|_T \in P^{k+1}(T), w \text{ continuous}\}. 
\end{eqnarray}
\end{comment}

\begin{comment}
$P^{k+1}(T)$ denotes a polynomial space of highest order $k$ on simplicial element $T$. We consider tetrahedral, hexahedral and  prismatic elements. For the standard continuous space $W_h$ and the tangential continuous space $V_h$ elements and shape functions are provided by Zaglmayr in \cite[p. \ 92ff.]{Zaglmayr:2006}. For the normal-normal continuous space elements and shape functions are described in \cite{PechsteinSchoeberl:11}.  In the open-source software package Netgen/Ngsolve (\href{https://ngsolve.org}{https://ngsolve.org} ) all elements are implemented. Additionally two-\-di\-men\-sional quadrilateral and triangular elements are provided. 
\end{comment}

The variational formulation of the elastic problem based on Reissner's principle is the following: find  $\uv$ satisfying $\uv_t = \mathbf{0}$ on $\Gamma_1$ and $\sigmat$ satisfying $\sigma_{nn} = 0$ on $\Gamma_2$ such that
\begin{align}
\int_\Omega \St :\sigmat : \delta \sigmat\, d\Omega - \langle \epst(\uv), \delta \sigmat\rangle -\langle \epst(\delta \uv), \sigmat \rangle &=\ \int_\Omega \fv \cdot \delta \uv\,d\Omega %+ %\int_{\Gamma_2} \tv_{nt} \delta \uv_t\, d\Gamma ,
\label{eq:TDNNS}
\end{align}
for all admissible virtual displacements $\delta \uv $ and virtual stresses $\delta\sigmat $ which satisfy the corresponding homogeneous essential boundary conditions. In (\ref{eq:TDNNS}) instead of the integral $\int_\Omega \epst(\uv) : \sigmat\, d\Omega$, the duality product $\langle \epst(\uv), \sigmat\rangle$ occurs. This is due to the fact that the displacement field $\uv$ is discontinuous, and therefore gaps between elements may arise. On finite element meshes, the duality product can be computed element wise by surface and volume integrals . Due to (additional) distributional parts of the strain, the surface integrals in (\ref{eq:defdiv1})-(\ref{eq:defdiv2}) do not vanish. Only if the stress field is normal-normal continuous, the duality product is well defined 

\begin{eqnarray}
\langle\epst(\uv), \sigmat \rangle &=& \sum_{T \in \mathcal{T}} \Big( \int_T \sigmat : \epst(\uv)\, d\Omega - \int_{\partial T} \sigmat_{nn}\cdot \uv_n\, d\Gamma \Big) \label{eq:defdiv1}\\
&=& \sum_{T \in \mathcal{T}} \Big( -\int_T \opdiv \sigmat \cdot \uv\, d\Omega + \int_{\partial T} \sigmat_{nt}\cdot \uv_t\, d\Gamma \Big)\ = \
-\langle \opdiv \sigmat, \uv\rangle. \label{eq:defdiv2}
\end{eqnarray}

Equations  (\ref{eq:defdiv1})-(\ref{eq:defdiv2}) can be shown by partial integration, using continuity of normal normal stress and tangential displacement and is shown in detail in \cite{PechsteinSchoeberl:11} and \cite{pechstein2018analysis}.
\par 

\subsection{Finite element formulation for piezoelastic problem }
To get equations for the (fully) coupled problem of piezoelasticity we  eliminate the dielectric displacement in the balance equations (\ref{eq:balance1})-(\ref{eq:balance2}) and use the constitutive equations in (\ref{eq:dtype1})-(\ref{eq:dtype2}) to get

\begin{eqnarray}
- \St^{E}: \sigmat + \dt^T\cdot \nabla \phi + \epst &= \mathbf{0}, \label{eq:dif1}\\
- \opdiv \sigmat &=\fv, \label{eq:dif2}\\
-\opdiv(\dt: \sigmat - \epsilont^{\sigma} \cdot\nabla \phi) &= 0. \label{eq:dif3}
\end{eqnarray}

 We  multiply (\ref{eq:dif1}) by a virtual stress $\delta \sigmat$, (\ref{eq:dif2}) by a virtual displacement $\delta \uv $ and (\ref{eq:dif3}) by a virtual electric potential $\delta \phi$ to get a variational formulation. The virtual quantities have to satisfy the corresponding homogeneous boundary conditions $\delta u_t=0$ on $\Gamma_1$, $\delta \sigma_{nn}=0$ on $\Gamma_2$ and $\delta\phi=0$ on $\Gamma_3$.  After integration %over the domain $\Omega$ 
using the identities in (\ref{eq:defdiv1}) - (\ref{eq:defdiv2}) we get from (\ref{eq:dif1}) and (\ref{eq:dif2})

\begin{eqnarray}
\int_\Omega (-\St^{E}: \sigmat + \dt^T\cdot \nabla \phi): \delta \sigmat\,d\Omega  + \langle\epst, \delta \sigmat\rangle &=& 0,\label{eq:varform1}\\
\label{eq:varform2} \langle \delta \epst, \sigmat\rangle &=& \int_\Omega \fv \cdot \delta \uv\, d\Omega.% + \int_{\Gamma_2} \tv_{nt} \cdot \delta \uv_t\, d\Gamma.
\end{eqnarray}

Integration by parts in (\ref{eq:dif3}), taking into account the homogeneous boundary conditions $\Gamma_3$ ($\delta\phi=0$)  and the natural  one on $\Gamma_4$ (no surface charges) leads to 

\begin{align}
\int_\Omega (\dt: \sigmat - \epsilont^{\sigma} \cdot\nabla \phi) \cdot \nabla\delta \phi\,d\Omega   = 0.
\end{align}

Summing up we finally arrive at 

\begin{align}
%-\int_\Omega (\St^{E} \sigmat - \dt^T\nabla \phi) : \delta \sigmat\, d\Omega 
%+ \langle \epst(\uv), \delta \sigmat\rangle + \langle \epst(\delta \uv), \sigmat \rangle +
%\int_\Omega (\dt \sigmat - \epsilont^{\sigma} \nabla \phi) \cdot \delta\nabla \phi\, d\Omega =\\
%\int_\Omega \Dv\cdot \delta \Ev\, d\Omega
-\int_\Omega (\St^{E}: \sigmat - \dt^T\cdot \nabla \phi) : \delta \sigmat\, d\Omega + \langle \epst(\uv), \delta \sigmat\rangle + \langle \epst(\delta \uv), \sigmat \rangle + \label{eq:V11}\\+\int_\Omega (\dt: \sigmat - \epsilont^{\sigma} \cdot\nabla \phi) \cdot \nabla \delta\phi\, d\Omega
=\label{eq:V12} \int_\Omega \fv \cdot \delta \uv\,d\Omega.% + \int_{\Gamma_2} \tv_{nt} \cdot \delta \uv_t\, d\Gamma . 
\end{align}

The performance of thin prismatic elements has been shown in \cite{PechsteinMeindlhumerHumer:2018}. They have been shown to be free from locking when using at least a order of $k=2$ for the electric potential and linear elements for the mechanical quantities (stress and displacement).

\section{Curvilinear elements} \label{sec:curvilinear_elements}

To implement the equations from the previous section in a computational code, we use the well known concept of reference elements. The reference element is of unit size and not distorted. Contrarily, elements in the finite element mesh are in general distorted, and may be curved in order to enhance geometry approximation. Integration of virtual work is then carried out after transformation to the reference element. Also shape functions are defined for the reference element.

In Section~\ref{sec_refelement}, we define the reference element for elements of prismatic and hexahedral shapes. Then we introduce shape functions for the stresses on these elements. Note that, for prismatic elements, we use less shape functions than proposed in the original paper \cite{PechsteinSchoeberl:12}, while shape functions for hexahedral elements are introduced for the first time. In Section~\ref{sec_transformation}, we describe in detail how quantities on the reference element are transformed to a distorted element in the mesh. Displacement shape functions are transformed by a covariant transformation, which preserves tangential degrees of freedom. From this transformation, we derive the correct transformation for the strain on a curvilinear element. Moreover, we present the correct transformation for the stress shape functions, as well as the transformation of the divergence of the stresses.
 %In the following, we describe in detail how tangential and normal-normal continuous shape functions are transformed to a possibly distorted, curvilinear element in such a way that the degrees of freedom are kept. Also the transformations of divergence and strain tensor are provided.

\subsection{Reference element and shape functions} \label{sec_refelement}

In the following, we will provide shape functions of arbitrary order for displacements, stresses and electric potential on the reference hexahedral and prism. As we use the tensor product structure of these elements, we first describe the reference segment and reference triangle. All quantities associated to reference elements are denoted by a hat. The shape functions are then transformed by according transformations described in Section~\ref{sec_transformation} to some possibly distorted element in the mesh.

For the electric potential, we use standard continuous hierarchical elements. For the displacement elements, we use tangential continuous elements well-known from electromagentics \cite{Nedelec:86}. %Degrees of freedom for potential and displacement elements for all element types are described in \cite{Monk:03}. 
All elements are implemented in the open-source software package Netgen/NGSolve \cite{netgen}. The exact basis function implemented there are described in detail in \cite{Zaglmayr:2006}. We keep close to the notation adopted in this latter reference. The stress basis functions were introduced for tetrahedra in \cite{PechsteinSchoeberl:11} and for prismatic elements in \cite{PechsteinSchoeberl:12}. In the following, we provide an adapted set of stress basis functions for prisms and stress basis functions for hexahedral elements.

\paragraph{The reference segment}
We use the unit reference segment $\hat T_{seg} = [0,1]$. For reference coordinate $\hat x$, we define
%\begin{itemize}
%\item 
the barycentric coordinates or ``hat basis functions'' 
\begin{align}
\hat \lambda_1(\hat x) &= 1-\hat x, & \hat \lambda_2(\hat x) &= \hat x,
\end{align}
%\item 
and the family of Legendre polynomials, $\hat q_i = \hat \ell_i(\hat x)$, where $i$ indicates the polynomial order.

%\end{itemize}
%\begin{comment}
\begin{figure}
\begin{center}
\psfrag{Ttrig}{$\hat T_{trig}$}
\psfrag{hatt}{$\hat T_{seg}$}
\psfrag{Tprism}{$\hat T_{prism}$}
\psfrag{Thex}{$\hat T_{hex}$}
\psfrag{lam1}{$\hat \lambda_1$}
\psfrag{lam2}{$\hat \lambda_2$}
\psfrag{x}{$\hat x$}
\psfrag{y}{$\hat y$}
\psfrag{z}{$\hat z$}
\psfrag{E1}{$\hat E_1$}
\psfrag{E2}{$\hat E_2$}
\psfrag{E3}{$\hat E_3$}
\psfrag{V1}{$\hat V_1$}
\psfrag{V2}{$\hat V_2$}
\psfrag{V3}{$\hat V_3$}
\psfrag{Fz1}{$\hat F^{\hat z}_1$}
\psfrag{Fz2}{$\hat F^{\hat z}_2$}
\psfrag{Fx1}{$\hat F^{\hat x}_1$}
\psfrag{Fx2}{$\hat F^{\hat x}_2$}
\psfrag{Fy1}{$\hat F^{\hat y}_1$}
\psfrag{Fy2}{$\hat F^{\hat y}_2$}
\psfrag{Fxy1}{$\hat F^{\hat x\hat y}_1$}
\psfrag{Fxy2}{$\hat F^{\hat x\hat y}_2$}
\psfrag{Fxy3}{$\hat F^{\hat x\hat y}_3$}
\includegraphics[width=0.7\textwidth]{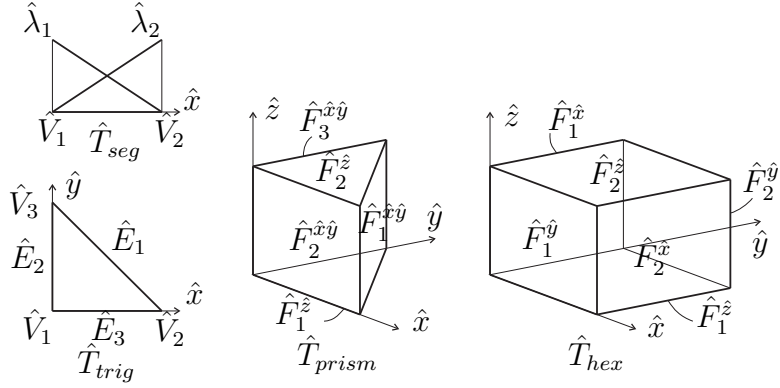}
\end{center}
\caption{The various unit elements.} \label{fig_refelements}
\end{figure}
%\end{comment}
\paragraph{The reference triangle}

The unit triangle shall be denoted by $\hat T_{trig} = \{(\hat x, \hat y): \hat x \geq0, \hat y \geq 0, \hat x + \hat y \leq 1\}$. For an enumeration of vertices and edges we refer to Figure~\ref{fig_refelements}. Also on the triangle, we define barycentric coordinates 
%\begin{itemize}
%\item 

\begin{align}
\hat \lambda_1(\hat x, \hat y) &= 1-\hat x-\hat y, & \hat \lambda_2(\hat x, \hat y) &= \hat x, & \hat \lambda_3(\hat x, \hat y) &= \hat y.
\end{align}
%\item 
We need a family of polynomials on the edge $E_\gamma$ between vertices $V_\alpha$ and $V_\beta$, that is extended to the whole triangle, e.g. the family of scaled Legendre polynomials
\begin{align}
\hat \ell^S_{\gamma,i}(\hat x, \hat y)  &= 
\hat \ell_i\left(\frac{\hat \lambda_\alpha - \hat \lambda_\beta}{\hat \lambda_\alpha  + \hat \lambda_\beta} \right) 
(\hat \lambda_\alpha  + \hat \lambda_\beta )^{i}
\end{align}
%\item 
A family of bivariate polynomials on the triangle of polynomial order $i+j$, can  e.g. be realized by
\begin{align}
\hat q_{ij}(\hat x, \hat y)  = \hat \ell^S_{3,i} \hat \ell_j(\hat \lambda_3 - \hat \lambda_1 - \hat \lambda_2).
\end{align}
Of course any other family of polynomials on the triangle may be used. The conditioning of the finite element matrices for high polynomial orders is affected by this choice, but this is not within the scope of the current contribution.
%\end{itemize}

On the triangle, we introduce three unit tensor fields, that have unit normal stress on one edge and zero normal stress on the other edges. Multiplying these unit tensors by (scalar) polynomials, we arrived at the desired hierarchical shape functions on the triangle, they can be found in \cite{PechsteinSchoeberl:11}. We set the unit tensors
\begin{align}
\hat \St_1 &= \left[ \begin{array}{cc} 0 & 1 \\ 1 & 0 \end{array} \right], &
\hat \St_2 &= \left[ \begin{array}{cc} 1 & -1/2 \\ -1/2 & 0 \end{array} \right], &
\hat \St_3 &= \left[ \begin{array}{cc} 0 & -1/2 \\ -1/2 & 1 \end{array} \right].
\end{align}

\paragraph{The reference prism}

The reference prism is constructed as a tensor product of the reference triangle (for $\hat x$ and $\hat y$ coordinates) and the reference segment (for $\hat z$ coordinate): $\hat T_{prism} = \hat T_{trig} \otimes \hat T_{seg} = \{(\hat x, \hat y, \hat z): (\hat x, \hat y) \in \hat T_{trig} \text{ and } \hat z \in \hat T_{seg}\}$. The reference prism is depicted in Figure~\ref{fig_refelements}.

We define the local reference shape functions for the element of order $p$.
On the prism, we have shape functions associated to each of the three quadrilateral faces $\hat F^{\hat x \hat y}_1$, $\hat F^{\hat x \hat y}_2$ and $\hat F^{\hat x \hat y}_3$, as well as shape functions associated to the two triangular faces $\hat F^{\hat z}_1$ and $\hat F^{\hat z}_2$, and shape functions that have zero normal stress and are associated with the element itself. The last class of shape functions is element-local and can be eliminated directly at assembly. We use
\begin{itemize}
\item for each quadrilateral face $\hat  F^{\hat x \hat y}_\gamma$, $(p+1)^2$ shape functions
\begin{align}
\hat \Nt^{\sigma,\hat F^{\hat x \hat y}_\gamma} _{ik} &= \left[ \begin{array}{cc}
\hat \ell_{\gamma,i}^S(\hat x, \hat y) \hat \ell_k(\hat z) \hat\St_\gamma & \nullv \\ \nullv & 0
\end{array}\right], & 
\begin{array}{l}0 \leq i,k \leq p, \\ \gamma = 1,2,3 .\end{array}
\end{align}
\item for each triangular face $\hat F^{\hat z}_\gamma$, $\frac12(p+1)(p+2)$ shape functions
\begin{align}
\hat \Nt^{\sigma,\hat F^{\hat z}_\gamma} _{ij} &= %\left[ \begin{array}{cc}
%\nullt & \nullv \\ \nullv & \hat q_{ij}(\hat x, \hat y)  \hat \lambda_{\gamma}(\hat z) 
%\end{array}\right]
\hat q_{ij}(\hat x, \hat y)  \hat \lambda_{\gamma}(\hat z)\ \ev_{\hat z} \otimes^s \ev_{\hat z}
, & \begin{array}{l}0 \leq i+j \leq p, \\ \gamma = 1,2 .\end{array}
\end{align}
\item three types of element-local shape functions with vanishing normal stress, each associated to one block of the symmetric three-by-three tensor,
\begin{align}
\hat \Nt^{\sigma,\hat x \hat y} _{\gamma ijk} &= \left[ \begin{array}{cc}
\hat q_{ij}(\hat x, \hat y)\hat \lambda_\gamma(\hat x, \hat y)  \hat \ell_j(\hat z) \hat\St_\gamma & \nullv \\ \nullv & 0
\end{array}\right], & 
\begin{array}{l}
0 \leq i+j \leq p-1,\\ 
0 \leq k \leq p+1,\\ 
\gamma = 1,2,3, \end{array}\\
\hat \Nt^{\sigma,\hat z}_{ijk} &= %\left[ \begin{array}{cc}
%\nullt & \nullv \\ \nullv & \hat q_{ij}(\hat x, \hat y)  \hat \lambda_{1}(\hat z) \hat \lambda_{2}(\hat z) \hat \ell_k(\hat z)
%\end{array}\right], & 
\hat q_{ij}(\hat x, \hat y)  \hat \lambda_{1}(\hat z) \hat \lambda_{2}(\hat z) \hat \ell_k(\hat z)\ \ev_{\hat z} \otimes^s \ev_{\hat z}, &
\begin{array}{l}
0 \leq i+j \leq p+1,\\ 
0 \leq k \leq p-1,\end{array}\\
\hat \Nt^{\sigma,\hat \xi \hat z} _{\hat \xi ijk} &= \hat q_{ij}(\hat x, \hat y) \hat \ell_k(\hat z)\ \ev_{\hat \xi} \otimes^s \ev_{\hat z}, & 
\begin{array}{l}
0 \leq i+j \leq p,\\ 
0 \leq k \leq p, \\
\hat \xi = \hat x, \hat y.\end{array}
\end{align}
\end{itemize}
Note that some components of the shape functions are of order $p+1$, however the element normal stresses are of order $p$.
\paragraph{The reference hexahedron}
The reference hexahedron is again constructed in tensor product style: $\hat T_{hex} = \hat T_{seg} \otimes \hat T_{seg} \otimes \hat T_{seg} $. %$= \{(\hat x, \hat y, \hat z): \hat x \in \hat T_{seg}, \hat y \in \hat T_{seg} \text{ and } \hat z \in \hat T_{seg}\}$. 
The reference hexahedron is depicted in Figure~\ref{fig_refelements}.

Again, we distinguish between shape functions that are associated with one of the quadrilateral faces, with zero normal stress on all other faces, and shape functions that have zero normal stress on all faces. Additionally, the shape functions are built in such a way that the  normal stress of the face-based shape functions corresponds to the normal stress distribution of face-based shape functions for prisms. This way, it is possible to use prismatic and hexahedral elements in the same mesh. Of course, also the shape functions for tetrahedral elements are implemented in NGSolve, which are not discussed here, match accordingly, such that hybrid meshes can be used.
We use the shape functions
\begin{itemize}
\item for each of the faces $\hat F^{\hat \xi}_\gamma$ with $\hat \xi \in \{\hat x, \hat y, \hat z\}$ and $\gamma = 1,2$
%\begin{align}
%\hat \Nt^{\sigma,F^{\hat x}_\gamma}_{jk} &= \left[ \begin{array}{ccc}
%\hat \ell_j(\hat y) \hat \ell_k(\hat z) \lambda_\gamma(\hat x) & 0 & 0 \\
%0 & 0 & 0 \\
%0 & 0 & 0
%\end{array} \right], & \begin{array}{l}0 \leq j,k \leq p, \\ \gamma = 1,2 .\end{array}
%\\
%\hat \Nt^{\sigma,F^{\hat y}_\gamma}_{ik} &= \left[ \begin{array}{ccc}
%0 & 0 & 0 \\
%0 & \hat \ell_i(\hat x) \hat \ell_k(\hat z) \lambda_\gamma(\hat y) & 0 \\
%0 & 0 & 0
%\end{array} \right], & \begin{array}{l}0 \leq j,k \leq p, \\ \gamma = 1,2 .\end{array}
%\\
%\hat \Nt^{\sigma,F^{\hat z}_\gamma}_{ij} &= \left[ \begin{array}{ccc}
%0 & 0 & 0 \\
%0 & 0 & 0 \\
%0 & 0 & \hat \ell_i(\hat x) \hat \ell_j(\hat y \lambda_\gamma(\hat z)
%\end{array} \right], & \begin{array}{l}0 \leq i,j \leq p, \\ \gamma = 1,2 .\end{array}
%\end{align}
\begin{align}
%\hat \Nt^{\sigma,F^{\hat x}_\gamma}_{jk} &= \hat \ell_j(\hat y) \hat \ell_k(\hat z) \lambda_\gamma(\hat x)\  \ev_{\hat x} \otimes^s \ev_{\hat x}, & \begin{array}{l}0 \leq j,k \leq p, \\ \gamma = 1,2,\end{array}
%\\
%\hat \Nt^{\sigma,F^{\hat y}_\gamma}_{ik} &= \hat \ell_i(\hat x) \hat \ell_k(\hat z) \lambda_\gamma(\hat y)\  \ev_{\hat y} \otimes^s \ev_{\hat y}, & \begin{array}{l}0 \leq j,k \leq p, \\ \gamma = 1,2,\end{array}
%\\
%\hat \Nt^{\sigma,F^{\hat z}_\gamma}_{ij} &= \hat \ell_i(\hat x) \hat \ell_j(\hat y) \lambda_\gamma(\hat z)\  \ev_{\hat z} \otimes^s \ev_{\hat z}, & \begin{array}{l}0 \leq i,j \leq p, \\ \gamma = 1,2 .\end{array}
\hat \Nt^{\sigma,F^{\hat \xi}_\gamma}_{ij} &= \hat \ell_i(\hat \eta) \hat \ell_j(\hat \zeta) \lambda_\gamma(\hat \xi)\  \ev_{\hat \xi} \otimes^s \ev_{\hat \xi}, & \begin{array}{l}0 \leq i,j \leq p, \\ \gamma = 1,2,\\ \{\hat \xi, \hat \eta, \hat \zeta\} = \{\hat x, \hat y, \hat z\}.\end{array}\end{align}
\item the element-local shape functions with vanishing normal stress
\begin{align}
%\hat \Nt^{\sigma,\hat x\hat x}_{ijk} &= \hat \ell_j(\hat y) \hat \ell_k(\hat z) \lambda_1(\hat x)\lambda_2(\hat x) \hat \ell_i(\hat x)\  \ev_{\hat x} \otimes^s \ev_{\hat x}, & \begin{array}{l}0 \leq j,k \leq p+1, \\ 0 \leq i \leq p-1,\end{array}
%\\
%\hat \Nt^{\sigma,\hat y \hat y }_{ijk} &= \hat \ell_i(\hat x) \hat \ell_k(\hat z) \lambda_1(\hat y)\lambda_2(\hat y) \hat \ell_j(\hat y)\  \ev_{\hat y} \otimes^s \ev_{\hat y}, & \begin{array}{l}0 \leq i,k \leq p+1, \\ 0 \leq j \leq p-1,\end{array}
%\\
%\hat \Nt^{\sigma,\hat z \hat z}_{ijk} &= \hat \ell_i(\hat x) \hat \ell_j(\hat y) \lambda_1(\hat z)\lambda_2(\hat z) \hat \ell_k(\hat z)\  \ev_{\hat z} \otimes^s \ev_{\hat z}, & \begin{array}{l}0 \leq i,j \leq p+1, \\ 0 \leq k \leq p-1,\end{array}
%\\
%\hat \Nt^{\sigma,\hat x\hat y}_{ijk} &= \hat \ell_i(\hat x)\hat \ell_j(\hat y) \hat \ell_k(\hat z) \  \ev_{\hat x} \otimes^s \ev_{\hat y}, & \begin{array}{l}0 \leq i,j \leq p, \\ 0 \leq k \leq p+1 ,\end{array}
%\\
%\hat \Nt^{\sigma,\hat y \hat z }_{ijk} &= \hat \ell_i(\hat x) \hat \ell_j(\hat y)\hat \ell_k(\hat z) \  \ev_{\hat y} \otimes^s \ev_{\hat z}, & \begin{array}{l}0 \leq j,k \leq p, \\ 0 \leq i \leq p+1,\end{array}
%\\
%\hat \Nt^{\sigma,\hat z \hat x}_{ijk} &= \hat \ell_i(\hat x) \hat \ell_j(\hat y) \hat \ell_k(\hat z)\  \ev_{\hat z} \otimes^s \ev_{\hat x}, & \begin{array}{l}0 \leq i,k \leq p, \\ 0 \leq j \leq p+1 .\end{array}
\hat \Nt^{\sigma,\hat \xi\hat \xi}_{ijk} &= \hat \ell_j(\hat \eta) \hat \ell_k(\hat \zeta) \lambda_1(\hat \xi)\lambda_2(\hat \xi) \hat \ell_i(\hat \xi)\  \ev_{\hat \xi} \otimes^s \ev_{\hat \xi}, & \begin{array}{l}0 \leq j,k \leq p+1, \\ 0 \leq i \leq p-1,\\ \{\hat \xi, \hat \eta, \hat \zeta\} = \{\hat x, \hat y, \hat z\},\end{array}\\
\hat \Nt^{\sigma,\hat \xi\hat \eta}_{ijk} &= \hat \ell_i(\hat \xi)\hat \ell_j(\hat \eta) \hat \ell_k(\hat \zeta) \  \ev_{\hat \xi} \otimes^s \ev_{\hat \eta}, & \begin{array}{l}0 \leq i,j \leq p, \\ 0 \leq k \leq p+1,\\ \{\hat \xi, \hat \eta, \hat \zeta\} = \{\hat x, \hat y, \hat z\}.\end{array}
\end{align}

\end{itemize}

\subsection{Transformations} \label{sec_transformation}

\begin{figure}[htp]
	\begin{center}
		\psfrag{Tdach}{$\widehat{T}$}
		\psfrag{xd1}{$1$}
		\psfrag{yd1}{$1$}
		\psfrag{xd}{$\widehat{x}_1$}
		\psfrag{yd}{$\widehat{x}_2$}
		\psfrag{x1}{$x_1$}
		\psfrag{y1}[r]{$x_2$}
		\psfrag{T}{$T$}
		\psfrag{phit}{$\Phi_T$,~$\Ft_T=\nabla\Phi_T$}
		\includegraphics[width=0.7\textwidth]{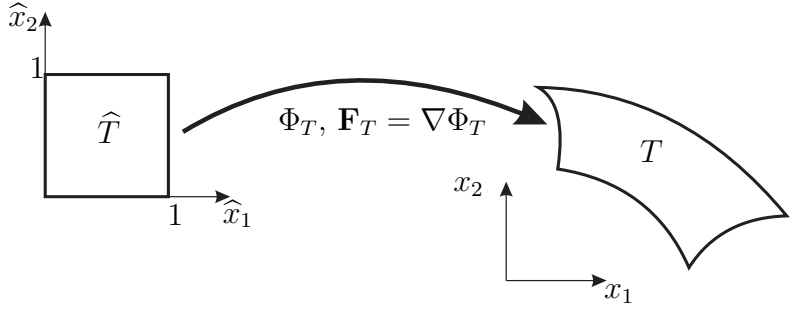}
	\end{center}
	\caption{Transformation from unit element to deformed element}
	\label{fig:tranform}
\end{figure}

Let now $\hat T$ denote the reference element of any type (triangular, quadrilateral, prismatic, tetrahedral or hexahedral), and let $T$ be a corresponding element in the finite element mesh. We assume $\Phiv_T(\hat{\xv})$ to be a smooth one to one mapping from $\hat{T}$ to $T$. A point in the reference element $\hat{\xv} \in \hat{T}$ is mapped by $\Phiv_T$ to some point $\xv\in T$. In Figure \ref{fig:tranform} this transformation is illustrated. Note, that in general  $\Phiv_T$ is nonlinear. The Jacobian of the transformation (similar to the deformation gradient tensor) is denoted by

\begin{equation}\label{eq:Ft}
	\Ft_T(\hat{x})=\nabla_{\hat{x}} \Phiv_T(\hat{\xv}). % =\left\{\frac{\partial \Phi_{T,i}}{\partial\hat{x}_j}(\hat{x})\right\}_{ij}.
\end{equation}

%\begin{comment}
%\end{comment}

%Here $\nabla_{\hat{x}}$ denotes the derivation with respect to the coordinates in the reference element $\hat{x}$ .
 The determinant of the Jacobian is denoted $J_T=\det(\Ft_T)$. Furthermore the  Hessian of the $i^{th}$ component of $\Phiv_T $, $\Ht^i_T$,  is given by

\begin{equation}\label{eq:hessian}
	(\Ht_T^i)_{jk}= \frac{\partial^2\Phiv_{T,i}}{\partial\hat{x}_j \partial\hat{x}_k}(\hat{\xv})  %,\qquad 1\leq i,j,k\leq d,
\end{equation}

% with $d$ the dimension of the problem. 
As already mentioned, the degrees of freedom of the method are $\uv_t$ and $\sigma_{nn}$. These degrees of freedom need to be preserved when shape functions are transferred to the curvilinear element. Assume $\hat{\Nt}^{\sigma}$ and $\hat{\Nt}^{u}$ to be one specific displacement and stress shape function, respectively. To calculate a finite element function, we transform these shape functions like

\begin{align}
	\Nt^u &=\Ft_T^{-T}\cdot\hat{\Nt}^u,\label{eq:transform_u}\\
	\Nt^{\sigma} &=\frac{1}{J_T^2}\Ft_T \cdot\hat{\Nt}^{\sigma} \cdot\Ft_T^T.\label{eq:transform_sigma}
\end{align}
%
%with $\Ft_T^{-T}=(\Ft_T^{-1})^T$. 
%Above, $\hat \uv$ and $\hat \sigmat$ denote some displacement respectively stress shape functions on the reference element, while $\uv$ and $\sigmat$ are the corresponding shape functions in the mesh.
The displacement and stress field are then weighted sums of these basis functions. \par 
 The first transformation for the displacement shape functions \eqref{eq:transform_u} is well known from finite elements for the electric field in Maxwell's equations, see e.g. \cite{Monk:03}. The linear strain $\epst$ of a finite element function is the the weighted sum of strains of such basis functions $\epst(\Nt^u)$. On curvilinear elements, this strain of a basis function has to be evaluated using the chain rule. In contrast to standard elements, it depends not only on the reference strain $\hat{\epst}(\hat{\Nt}^u)=\frac{1}{2}\left(\nabla_{\hat{\xv}}\hat{\Nt}^u+\nabla_{\hat{\xv}}\hat{\Nt}^{u,T}\right)$ but also on the shape function $\hat{\Nt}^u$ %As a consequence of (\ref{eq:transform_u}) and the chain rule, the strain  $\epst(\uv)$ is given by

\begin{equation}\label{eq:strain_u}
	\epst(\Nt^u)=\Ft_T^{-1}\cdot\hat{\epst}(\hat{\Nt}^u)\cdot\Ft_T^{-T} + \sum_{i=1}^{d}\hat{\Nt}^u_i\Ft_T^{-1}\cdot(\Ht^i_T)^{-1}\cdot\Ft_T^{-T}.
\end{equation}

The second part of (\ref{eq:strain_u}) only vanishes if $\Ft_T$ is constant, which is the case only for affine linear elements.\par
The transformation for stress shape functions (\ref{eq:transform_sigma}) can be found in \cite{PechsteinSchoeberl:11}. Note that compared to the Piola transformation, we have a factor of $\frac{1}{J_T^2}$ instead of $\frac{1}{J_T}$.
 For the (mechanical) balance equation (\ref{eq:balance1}) the divergence of the stress tensor has to be derived. Again, it consists of a weighted sum of $\opdiv_{\xv}(\Nt^{\sigma})$. From basic calculus for the standard Piola transformation, we deduce 

\begin{equation}
	 \opdiv_{\xv}(\Nt^{\sigma})=\opdiv_{\xv}\left(\frac{1}{J_T^2}\Ft_T\cdot\hat{\Nt}^{\sigma}\cdot\Ft_T^T\right)=\frac{1}{J_T}\opdiv_{\hat{\xv}}\left(\tilde{\Ft}_T\cdot\hat{\Nt}^{\sigma}\right),
 \end{equation}
 where $\tilde\Ft_T:=\frac{1}{J_T}\Ft_T$. By application of the product rule we get the $i^{th}$ component of the divergence of the stress tensor as
\begin{equation}
	\label{eq:div_sigma}
	\left(\opdiv_{\xv}\Nt^{\sigma}\right)_i=\frac{1}{J_T^2}\Ft_{T,ij}\left(\opdiv_{\hat{\xv}}\hat{\Nt}^{\sigma}\right)_j+\frac{1}{J_T}\frac{\partial\tilde{\Ft}_{T,ik}}{\partial\hat{x}_j}\hat{\Nt}^{\sigma}_{kj}.
\end{equation}
 
 \begin{comment}
 From the Piola transformation for the divergence of a tensor $\taut=\frac{1}{J_T}\hat{\taut}~ \Ft_T$ in the actual configuration the identity

\begin{equation}\label{eq:piolatransfomation}
	\opdiv_x(\taut)=\opdiv_x \left(\frac{1}{J_T}\cdot\hat{\taut}~ \Ft_T^T\right)=\frac{1}{J_T}\opdiv_{\hat{x}}(\hat{\taut})
\end{equation}

is known. %Analogous to $\nabla_{\hat{x}}$, the divergence with respect to $\hat{x}$ is denoted  $\opdiv_{\hat{x}}$.
 We introduce $\tilde{\Ft}_T=\frac{1}{J_T}\Ft_T$ and use (\ref{eq:piolatransfomation}) and the transformation (\ref{eq:transform_sigma})  get the divergence of the stress tensor 

\begin{eqnarray}
	\opdiv_x(\sigmat) &=& \opdiv_x \left( \frac{1}{J_T}\tilde{\Ft}_T\cdot\hat{\sigmat} \cdot\Ft_T^T \right)  \\
	&=& \frac{1}{J_T} \opdiv_{\hat{x}}\left(\tilde{\Ft}_T \cdot\hat{\sigmat} \right)\label{eq:div_sigma2}
\end{eqnarray}

in the reference system $\hat{x}$. We further treat (\ref{eq:div_sigma2}) by application of the product rule. We get the $i^{th}$ component of the divergence of the stress tensor , using the summation convention,

\begin{equation}
(\opdiv_x(\sigmat))_i = \frac{1}{J_T^2}\Ft_{T,ij}(\opdiv_{\hat{x}}\hat{\sigmat})_j + \frac{1}{J_T}\frac{\partial\tilde{\Ft}_{T,ik}}{\partial \hat{x}_j}\hat{\sigmat}_{kj}.
\end{equation}
\end{comment}
Assuming that the shape functions $\hat{\Nt}^{\sigma}$ as well as their divergence $\opdiv_{\hat{\xv}}\hat{\Nt}^{\sigma}$ are known analytically on the reference element, formula (\ref{eq:div_sigma}) allows to evaluate $\opdiv_{\xv}\Nt^{\sigma}$ for a curvilinear mesh element. 

\section{Eigenvalue problem} \label{sec:eigenvalues}
A wide range of technical applications involving piezoelectric structures is operated in the regime of (mechanical) vibrations e.g. noise and vibration control and damping, shape control or energy harvesting. To design and study these applications the computation of their  eigenfrequencies and -forms may be required. % They may be required for the design of feedback control systems \cite{moheimani2006piezoelectric}.
The finite element discretization leads to a system with a huge number of eigenvalues of which typically only a few of the smallest are of interest. In our contribution we propose to use the inverse iteration procedure which allows the effective computation of a certain number of eigenvalues. 
%\red{todo, Einleitung dazu} \par

\subsection{Inverse iteration}

The inverse iteration procedure \cite{schwarz2006NuMa,KNYAZEV2003,NEYMEYR200161} is designed to find the smallest eigenvalues and -vectors of the generalized eigenvalue problem
\begin{equation}\label{eq:ev_problem}
	\Am\cdot\qv_k=\lamv_k\Mm\cdot\qv_k
\end{equation}
with symmetric positive definite stiffness matrix $\Am$ and  mass matrix $\Mm$. %Here, $\Am$ and $\Mm$ are assumed to be symmetric positive definite matrices. 

%By the inverse iteration procedure, one may find the smallest eigenvalues and -vectors of the generalized eigenvalue problem
%For the theoretical part we treat the standard problem by setting $\Mm=\Imm$, with $\Imm$ the identity matrix, without restrictions of the original problem (\ref{eq:ev_problem}).  
%well known 
For a given $\lambda_k^n$ and $\qv_k^n$, the next iterate $\qv_k^{n+1}$ is computed by solving a linear system

\begin{equation}\label{eq:inv_iter}
	\Am \cdot\qv^{n+1}_k=\lamv_k^n\Mm\cdot\qv^n_k.
\end{equation}

%maps the actual eigenvector iterate $\xv^n_k$ to the next iterate $\xv^{n+1}_k$ by solving a system of linear equations. 
An approximation of the eigenvalue is given by the Rayleigh quotient
\begin{equation}\label{eq:rayleight_quotient}
\lamv(\qv)=\frac{\qv^T\cdot\Am\cdot\qv}{\qv^T\cdot\Mm\cdot\qv}.
\end{equation}
 For a suitable start vector $\qv^0$, which is not perpendicular to the smallest eigenvector, the inverse iteration procedure converges to the smallest eigenvalue $\lamv_1$ and the corresponding Eigenvector $\qv_1$.

\subsection{Structure of assembled system}
%\red{3x3-Matrix Struktur von $\mathbf{K}$}

The TDNNS (mixed) method for piezoelasticity with the variational formulation (\ref{eq:V11}-\ref{eq:V12}) leads to a certain structure of the global assembled system. The global stiffness matrix $\Am_{global}$ is indefinite, the global mass matrix $\Mm_{global}$ is positive semidefinite. The eigenvalue problem  is of the form

\begin{comment}
\begin{eqnarray*}
\drawmatrix[size=1.8]{ \mathbf{0} }\drawmatrix[size=1.8, width=1.3]{\mathbf{B}^T} \drawmatrix[size=1.8,width=.9]{\mathbf{0}} \hspace{15pt} \drawmatrix[size=1.8,width=1]{\uv}\hspace{15pt}  \drawmatrix[size=1.8,width=1]{\lambda \mathbf{M}}\\
\drawmatrix[size=1.3,width=1.8]{ \mathbf{B} }\drawmatrix[size=1.3]{-\mathbf{C}} \drawmatrix[size=1.3,width=.9]{\mathbf{D}^T} \times \drawmatrix[size=1.3,width=1]{\mathbf{\sigma}}=\drawmatrix[size=1.3,width=1]{\mathbf{0}}  \\
\drawmatrix[size=.9,width=1.8]{ \mathbf{B} }\drawmatrix[size=.9,width=1.3]{\mathbf{D}} \drawmatrix[size=.9]{\mathbf-{E}} \hspace{15pt}\drawmatrix[size=0.9,width=1]{\mathbf{\Phi}}\hspace{15pt}  \drawmatrix[size=0.9,width=1]{\mathbf{0}}
\end{eqnarray*}
\end{comment}

\begin{eqnarray}
\underbrace{\left[\begin{BMAT}(e)[1pt,3cm,3cm]{c:c:c}{c:c:c}
\nullt&\Bm^T&\nullt\\
\Bm&-\Cm&\Dm^T\\
\nullt&\Dm&-\Em
\end{BMAT}\right] \cdot}_{\Am_{global}}
\underbrace{\left[\begin{BMAT}(e)[1pt,.5cm,3cm]{c}{c:c:c} 
\qv_{\uv}\\ \qv_{\sigmat }\\ \qv_{\phim}\end{BMAT}\right]}_{\qv} = 
\lambda\underbrace{\left[\begin{BMAT}[1pt,3cm,3cm]{c:c:c}{c:c:c} 
 \Mm &\nullt& \nullt \\  \nullt & \nullt &\nullt \\ \nullt&\nullt & \nullt
\end{BMAT}\right]}_{\Mm_{global}}\cdot\underbrace{\left[\begin{BMAT}(e)[1pt,.5cm,3cm]{c}{c:c:c} 
\qv_{\uv}\\ \qv_{\sigmat }\\ \qv_{\phim}\end{BMAT}\right]}_{\qv}.
\label{eq:structur}
\end{eqnarray}
The global vector of degrees of freedom is denoted by $\qv$, where $\qv_{\uv}$, $\qv_{\sigmat }$ and $\qv_{\phim}$ denote the degrees of freedom of displacement, stress and electric potential, respectively.
The symmetric positive definite matrices $\Cm$ and $\Em$ correspond to the mechanical compliance and electrical permittivity of the constitutive equations. In the variational formulation these parts are represented by $\int_{\Omega}\sigmat:\St^E:\delta\sigmat d\Omega$ and $\int_{\Omega}\nabla\phi\cdot\epsilont\cdot\delta\nabla\phi d\Omega$. The matrix $\Dm$ corresponds the electromechanical coupling terms $\int_{\Omega}\left( \dt: \sigmat\cdot\nabla \delta\phi\right)d\Omega$, while $\Bm$ represents $\langle \epst(\uv), \delta \sigmat\rangle$, defined in (\ref{eq:defdiv1})-(\ref{eq:defdiv2}). Note that in (\ref{eq:structur}) the right hand side is $\mathbf{0}$ for the stresses and the electric potential. Therefore the global mass matrix is positive semidefinite. \par 
%In order to use the inverse iteration procedure which requires positive definite stiffness and mass matrices. Therefore the electric potential and the stresses in (\ref{eq:structur}) have to be eliminated. 
In the following, we motivate why the inverse iteration procedure can be used although neither the mass nor the stiffness matrix are positive definite. We emphasize that these considerations are purely theoretical, while the procedure is then directly applied to the indefinite or semidefinite matrices in implementations.  
From the last line of (\ref{eq:structur}) the potential degrees of freedom can be computed as

\begin{equation}\label{eq:phi_red}
	\qv_{\phim}=\Em^{-1}\cdot\Dm\cdot\qv_{\sigmat }.
\end{equation}
Inverting (\ref{eq:phi_red}) into the second line of (\ref{eq:structur}) we get for the stresses
\begin{equation}\label{eq:sigma_red}
	\Bm\cdot\qv_{\uv}=\left(\Cm-\Dm^T\cdot\Em^{-1}\cdot\Dm\right)\qv_{\sigmat }=\bar{\Cm}\cdot\qv_{\sigmat }.
\end{equation}
Finally applying (\ref{eq:sigma_red}) to the fist line of (\ref{eq:structur}) we get the eigenvalue problem
\begin{equation}\label{eq:ev_problem_struct}
	\Bm^T\cdot\bar{\Cm}^{-1}\cdot\Bm\cdot\qv_{\uv}=\lambda\Mm\cdot\qv_{\uv}.	
\end{equation}
For typical technical piezoelectric materials the matrix $\bar{\Cm}$ can be assumed to be positive definite as it represents the compliance at constant dielectric displacement. Therefore, (\ref{eq:sigma_red}) is formally equivalent to a generalized eigenvalue problem with symmetric positive definite matrices, and can be treated by the inverse iteration procedure.% the eigenvalue problem (\ref{eq:ev_problem_struct}) can be treated by the inverse iteration procedure. 
\par 
Note that the elimination of the degrees of freedom of stress an potential in (\ref{eq:phi_red})-(\ref{eq:ev_problem_struct}) is only of theoretical interest, to show that the eigenvalue problem can be solved by the inverse iteration procedure. In the practical implementation the global stiffness and mass matrix  $\Am_{global}$ and $\Mm_{global}$ are used. %Equations (\ref{eq:phi_red}) and (\ref{eq:sigma_red}) are solved within the iteration procedure. 
An exemplary implementation of the inverse iteration procedure \textit{Netgen/NGSolve}  can be found in the tutorial section in \cite{netgen}.

\section{Numerical Results} \label{sec:num_results}
In the sequel we study two numerical examples. First we show the discretization of a circular piezoelectric patch actor applied to a square aluminium plate. We evaluate the displacement of the plate and the stress field for an applied static voltage as well as the eigenfrequencies and eigenforms of the assembly. Our second example is a radially polarized piezoelectric half cylinder, studied in \cite{ZOUARI:2010}. We calculate (static) displacements and present a convergence study. %We show how the transformation for shape functions can be used in order to transform material laws. Finally we show how transformation laws can be used in order to consider (possibly large) predeformations and eigenstrains. 
For both examples we compare our results to results gained in the commercial software tool \textit{ABAQUS} 6.14 \cite{abaqus20146}.

\subsection{Circular patch actor}
\begin{figure}[htp]

\centering
\begin{subfigure}[b]{0.45\textwidth}
	\centering
		\psfrag{x}{$x$}
		\psfrag{y}{$y$}
		\psfrag{z}{$z$}
		\psfrag{voltage}{\footnotesize $\Delta\phi=100$ V}
		\psfrag{back}{\footnotesize fixed}
		\psfrag{def}{\footnotesize deformation}
		\psfrag{piezo}{\begin{tabular}{@{}l@{}}
				\footnotesize piezoelectric patch actor
		\end{tabular}}
		\psfrag{aluminium}{\footnotesize{\begin{tabular}{@{}l@{}}
				 aluminium \\
				 plate
		\end{tabular}}}
		
		\includegraphics[width=0.95\textwidth]{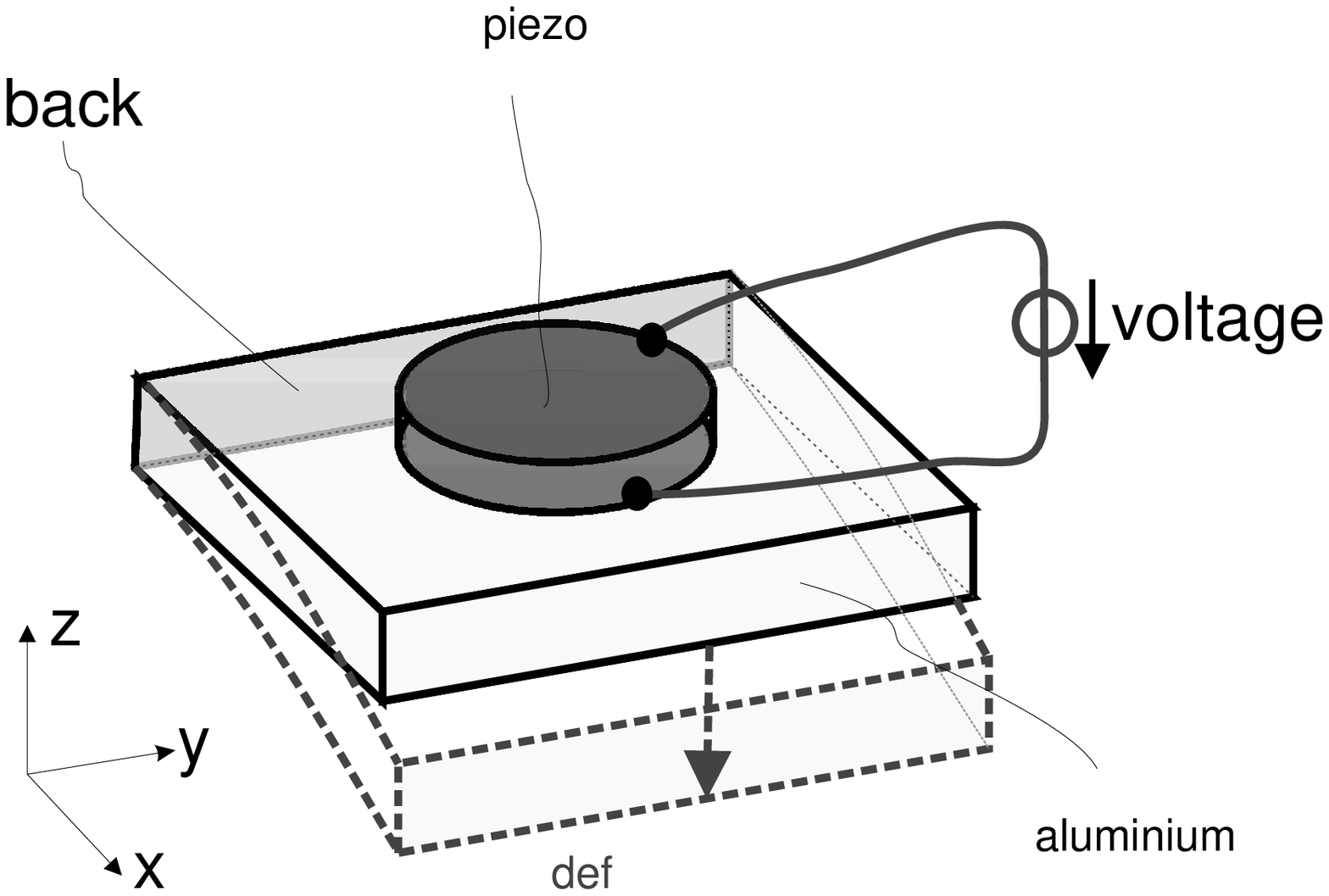}

	\caption{\centering Sketch of circular patch on quadratic aluminium plate}
	\label{fig:scetch_circ_patch}
	
\end{subfigure}\hspace{1eM}
\begin{subfigure}[b]{0.45\textwidth}
	
	\centering
	\includegraphics[width=0.85\textwidth]{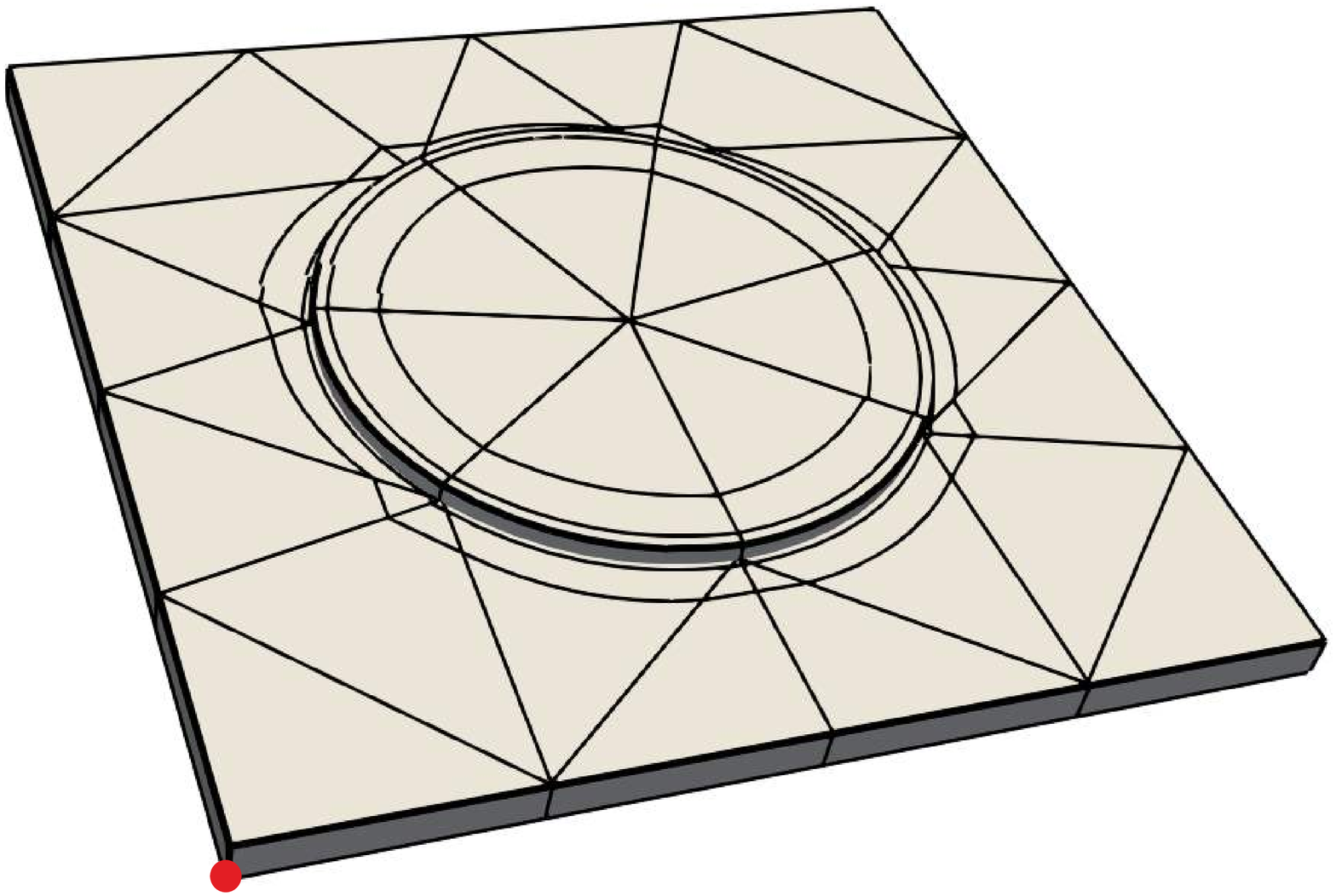}
	
	\caption{\centering  Mesh for calculation in \textit{Netgen/NGSolve}  ($\approx 12k$ DOF) }
	\label{fig:mesh_circ_patch}
\end{subfigure}

\caption{Sketch and mesh for circular piezoelectric patch actor}
\label{fig:circPatch_sketch_and_mesh}
\end{figure}

Our first numerical example is a circular piezoelectric patch  applied on a square plate. A schematic sketch is shown in Figure \ref{fig:circPatch_sketch_and_mesh}. The piezoelectric patch has a diameter of $15\operatorname{mm}$ and a thickness of $0.5\operatorname{mm}$. The patch material is considered to be PTZ-H5, polarized in thickness direction. The material properties are taken from \cite{ZOUARI:2010} and summarized in Table \ref{tab:mat_piezp}. The  square aluminium plate (Young's modulus $65\operatorname{GPa}$, Poisson ratio $0.3$, density $2.7 \times 10^{-9} \operatorname{kg/mm^3}$) has a length of $25 \operatorname{mm}$ and a thickness of $1\operatorname{mm}$. One side of the plate is clamped.  A potential difference  $\Delta\phi=100\operatorname{V}$ is applied to the electrodes of the piezoelectric patch, which are located at the top and bottom of the patch actor.

%\begin{comment}
%\end{comment}

%\begin{comment}
\begin{table}[htp]
	\centering
	\caption{Elastic ($[C]=\operatorname{GPa}$), piezoelectric ($[e]=\operatorname{C/mm^2} $) and dielectric ($[\epsilon]=\operatorname{C/Vm}$) properties and density ($[\rho]=\operatorname{kg/mm^3}$) of PZT-5H}
	\begin{tabular}{|l  l|}
		\hline
		$C_{11}^E=127.205 $ & $e_{31}=-6.62\times 10^{-3}$  \\ 
		$C_{12}^E=80.212 $ & $e_{33}=23.24\times 10^{-3}$ \\ 
		$C_{13}^E=84.670 $ & $e_{15}=17.03\times 10^{-3}$ \\ 
		$C_{33}^E=117.436 $ & $\epsilon_{11}=-6.62\times 10^{-9}$ \\ 
		$C_{66}^E=23.474 $ & $\epsilon_{33}=-6.62\times 10^{-9}$ \\ 
		$C_{44}^E=22.988 $ &  $\rho=7.5\times10^{-9}$\\ \hline
	\end{tabular}	
	\label{tab:mat_piezp}
\end{table}
%\end{comment}

 For our reference calculations in \textit{ABAQUS}  we use a mesh with an element size $0.25\operatorname{mm}$, consisting of four elements across the thickness of the aluminium and two elements across the thickness of the patch. We use \textit{ABAQUS} 6.14 C3D20R elements (20 nodes quadratic brick elements, reduced integration)  for the linear elastic aluminium plate and C3D20RE elements (20 node quadratic piezoelectric brick, reduced integration) for the piezoelectric patch. Our \textit{ABAQUS} model requires approximately $830k$ degrees of freedom.
 \par
 The mesh used for the calculation in \textit{Netgen/NGSolve} %\cite{netgen}
  is shown in Figure \ref{fig:mesh_circ_patch}. It consists of a very coarse prismatic mesh using only one element in thickness direction, for both the plate and the patch. % The patch surface is meshed with six elements prismatic elements.
   On the circumference of the patch the mesh is refined by two ring layers of hexahedral elements, which have an aspect ration of $\approx 25:2:1$. We use shape functions of order $p=2$ for displacements and stresses and order $p_{\phi}=3$ for the electric potential.  Summing up, we use $12439$ degrees of freedom for calculations in \textit{Netgen/NGSolve}.
   
   \par Note that in the TDNNS-method, contrary to elements using nodal degrees of freedom, a reduction of the thickness of the patch (or the plate) has no effect on the  number of used degrees of freedom. Actually we have chosen a rather thick patch-actor, in order to get accurate reference solutions at reasonable computational costs. 
   
\par
\begin{figure}[htp]
	\centering
	\begin{subfigure}[t]{0.4\textwidth}
		\includegraphics[width=0.8\textwidth]{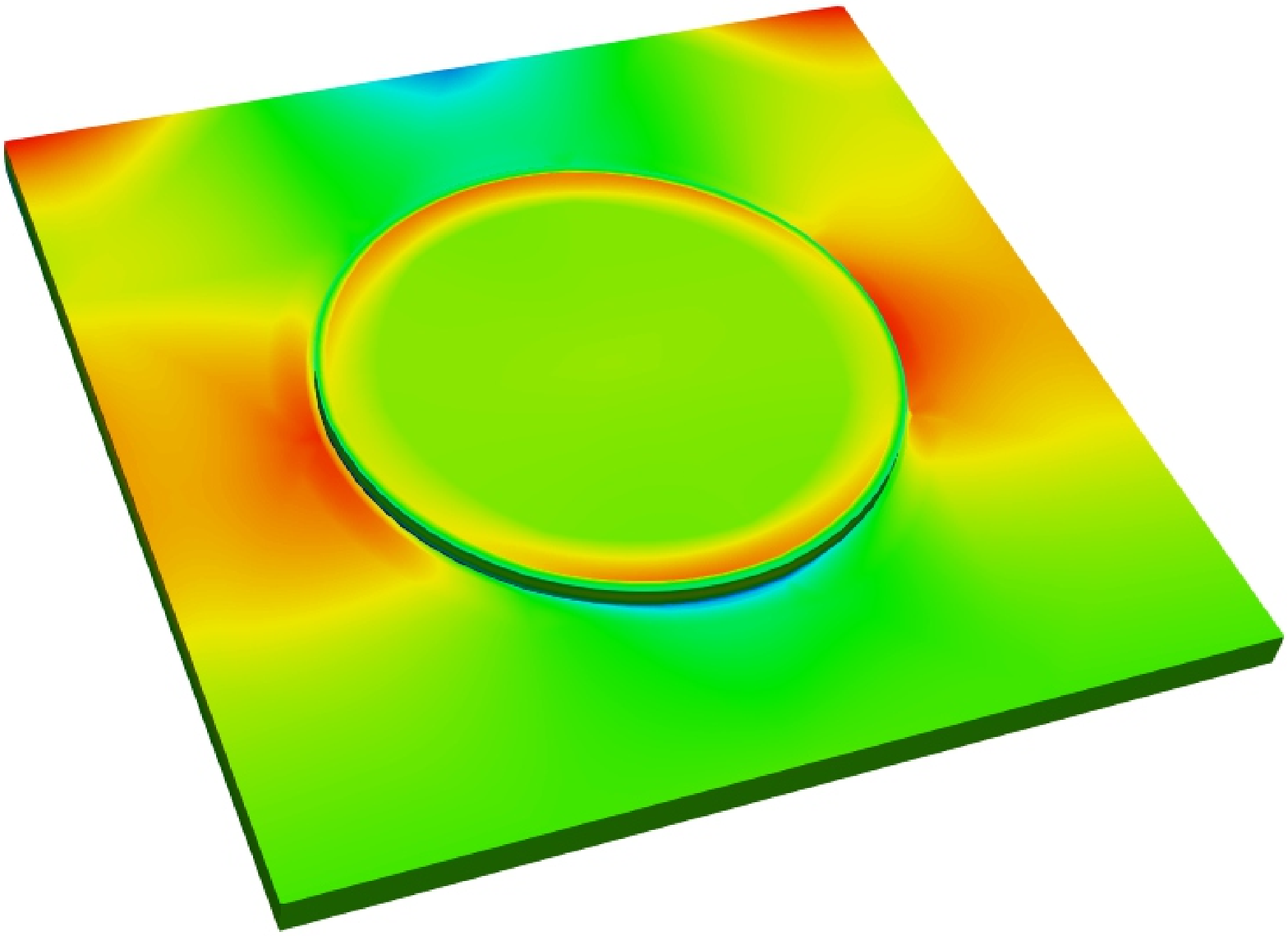}
		\subcaption{Solution \textit{Netgen/NGsolve} \\ ($\approx12k$ DOF)}
		\label{fig:comp_stress_netgen}
	\end{subfigure}
	\begin{subfigure}[t]{0.4\textwidth}
		\includegraphics[width=.8\textwidth]{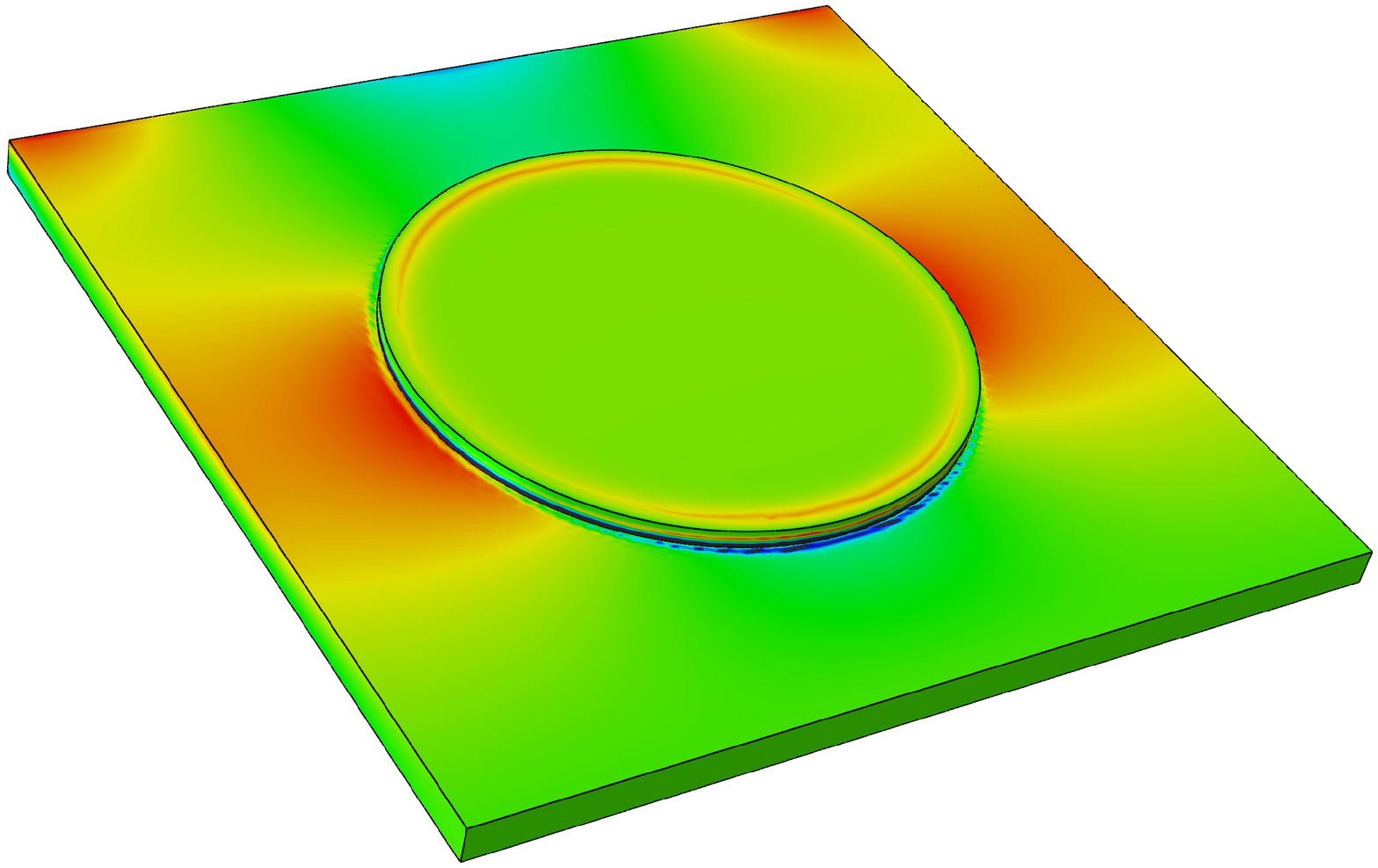}
		\subcaption{ Solution \textit{ABAQUS} \\ ($\approx 830k$ DOF)}
		\label{fig:comp_stress_abaqus}
	\end{subfigure}
	\begin{subfigure}[t]{0.15\textwidth}
		\psfrag{stress}{$S_{11}$}
		\psfrag{-2.0}{\footnotesize$-2.0$}
		\psfrag{-1.0}{\footnotesize$-1.0$}
		\psfrag{0.0}{\footnotesize$0$}
		\psfrag{1.0}{\footnotesize$1.0$}
		\psfrag{1.5}{\footnotesize$1.5$}
		\includegraphics[width=0.9\textwidth]{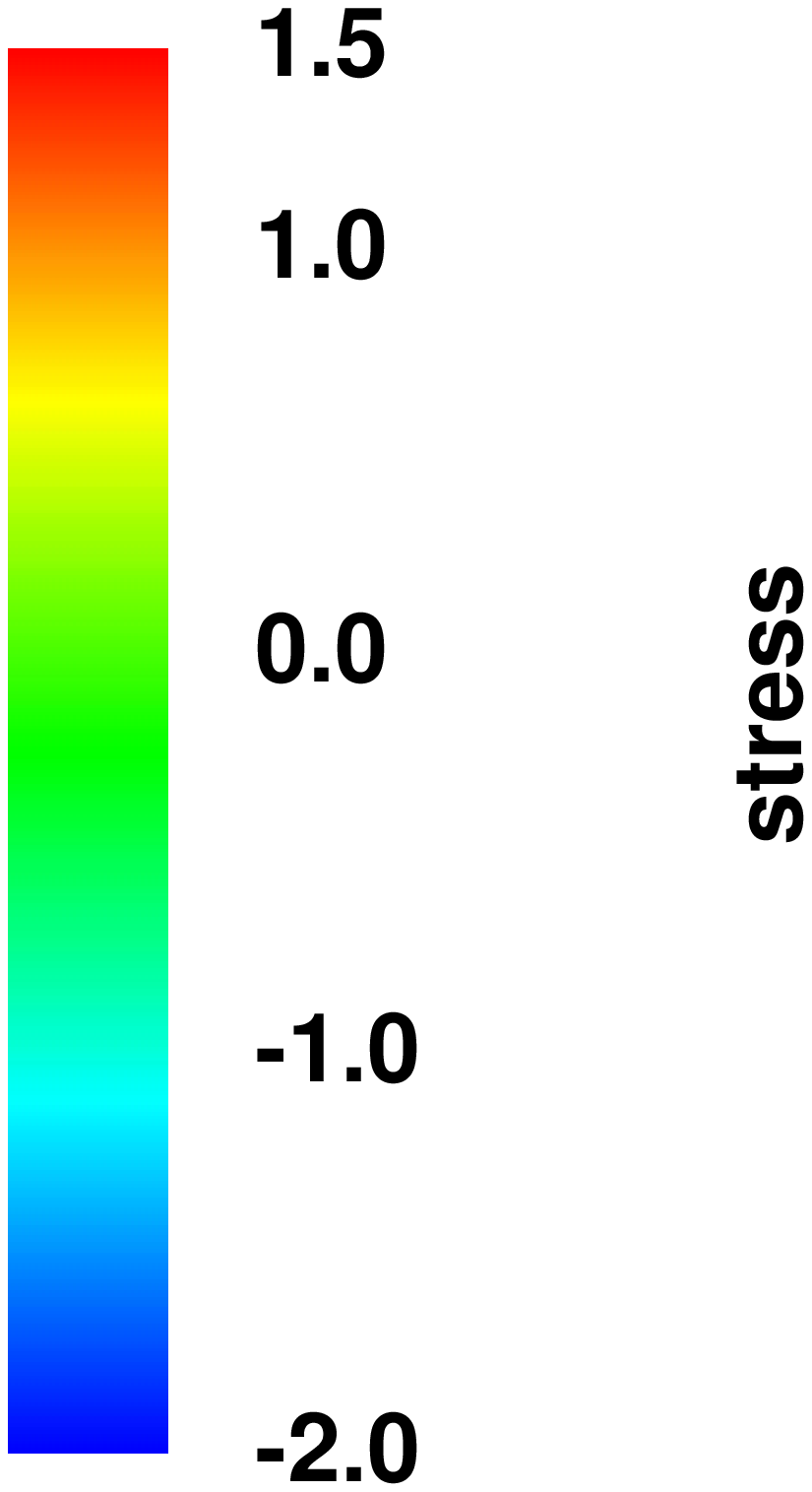}
	\end{subfigure}
	\caption{Contour plot of stress $S_{11}$ in \textit{ABAQUS} and \textit{Netgen/NGSolve}}
	\label{fig:compare_stess_circ_patch}
\end{figure}

 We compare resulting stress fields. In Figure \ref{fig:compare_stess_circ_patch} show  contour plots of the stress component $S_{11}$. In Figure \ref{fig:compare_stess_circ_patch} the stress fields  calculated in \textit{Netgen/NGSolve} and reference solution are shown.
 We as well evaluate the resulting vertical displacements of the corner point marked in Figure \ref{fig:circPatch_sketch_and_mesh}(\subref{fig:mesh_circ_patch}). While the reference displacement is $-6.07325\operatorname{\mu m}$, our result is $-6.03698\operatorname{\mu m}$. The relative difference of the results is $-0.597\%$.

%\begin{comment}
%\end{comment}

\begin{figure}[htp]
	\centering
	\begin{subfigure}[t]{0.3\textwidth}
		\includegraphics[width=0.8\textwidth]{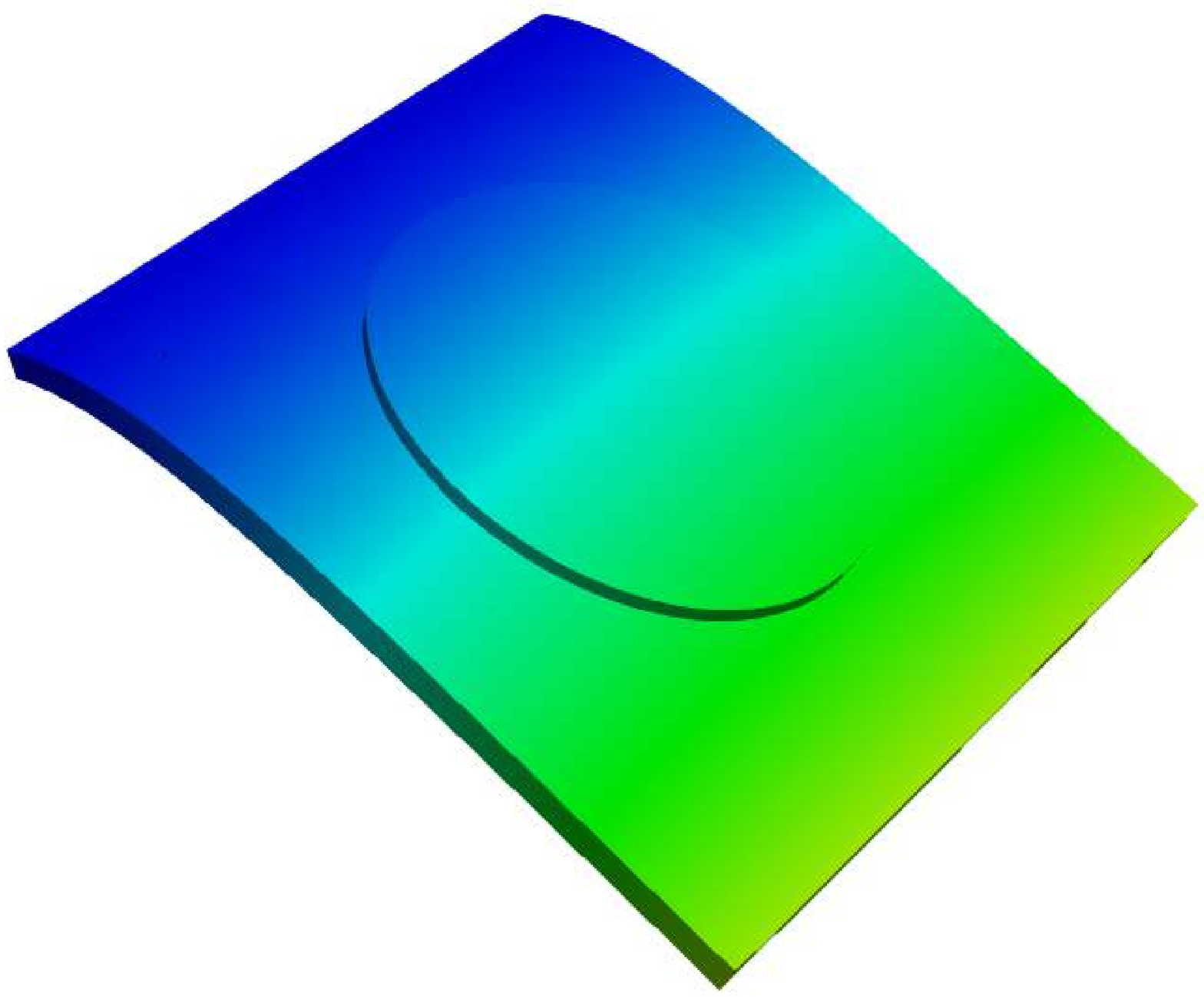}
%		\subcaption{First Eigenform}
%		\label{fig:circPatch_EF1}
	\end{subfigure}
	\begin{subfigure}[t]{0.3\textwidth}
		\includegraphics[width=0.8\textwidth]{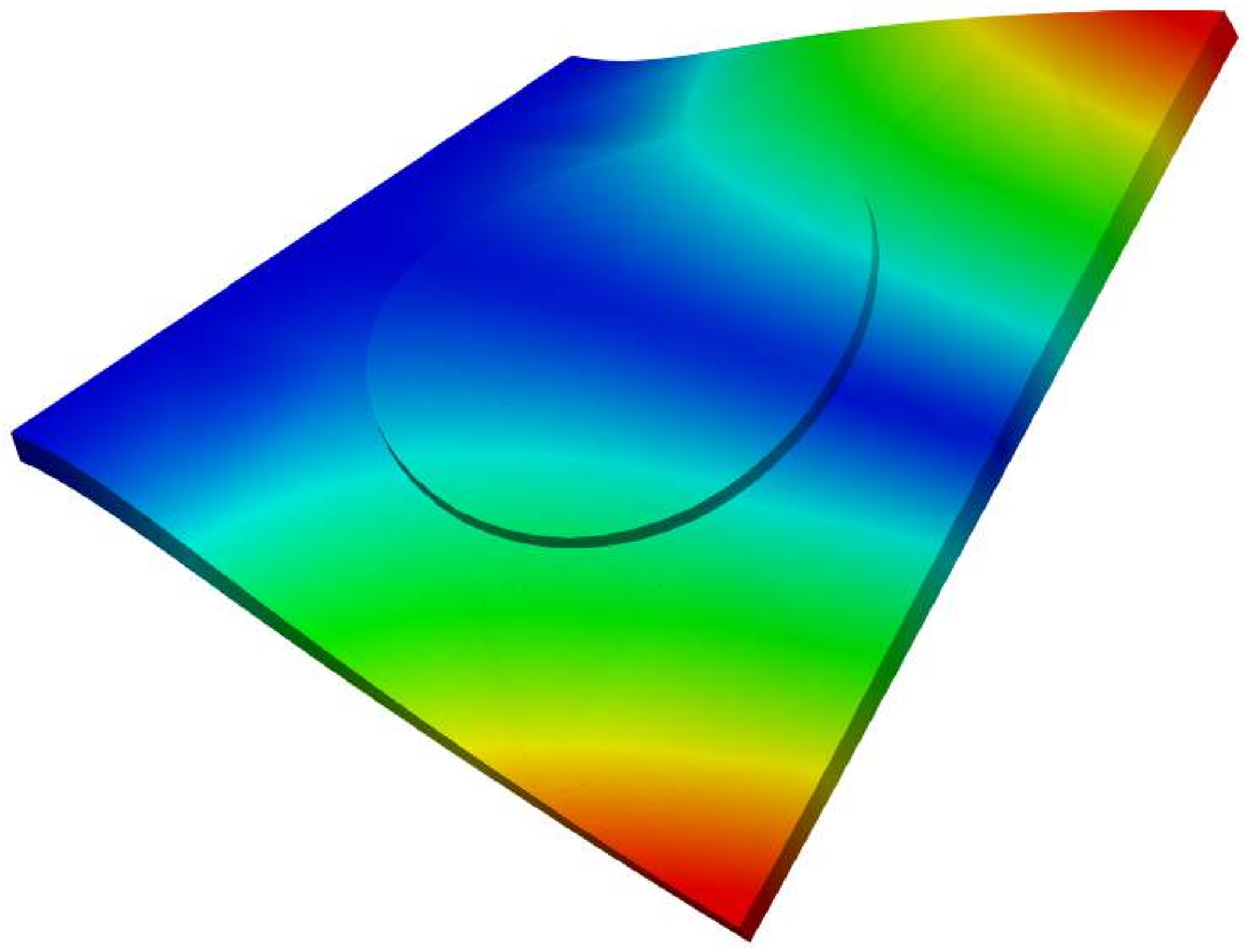}
%		\subcaption{Second Eigenform}
%		\label{fig:circPatch_EF2}
	\end{subfigure}
	\begin{subfigure}[t]{0.3\textwidth}
		\includegraphics[width=0.8\textwidth]{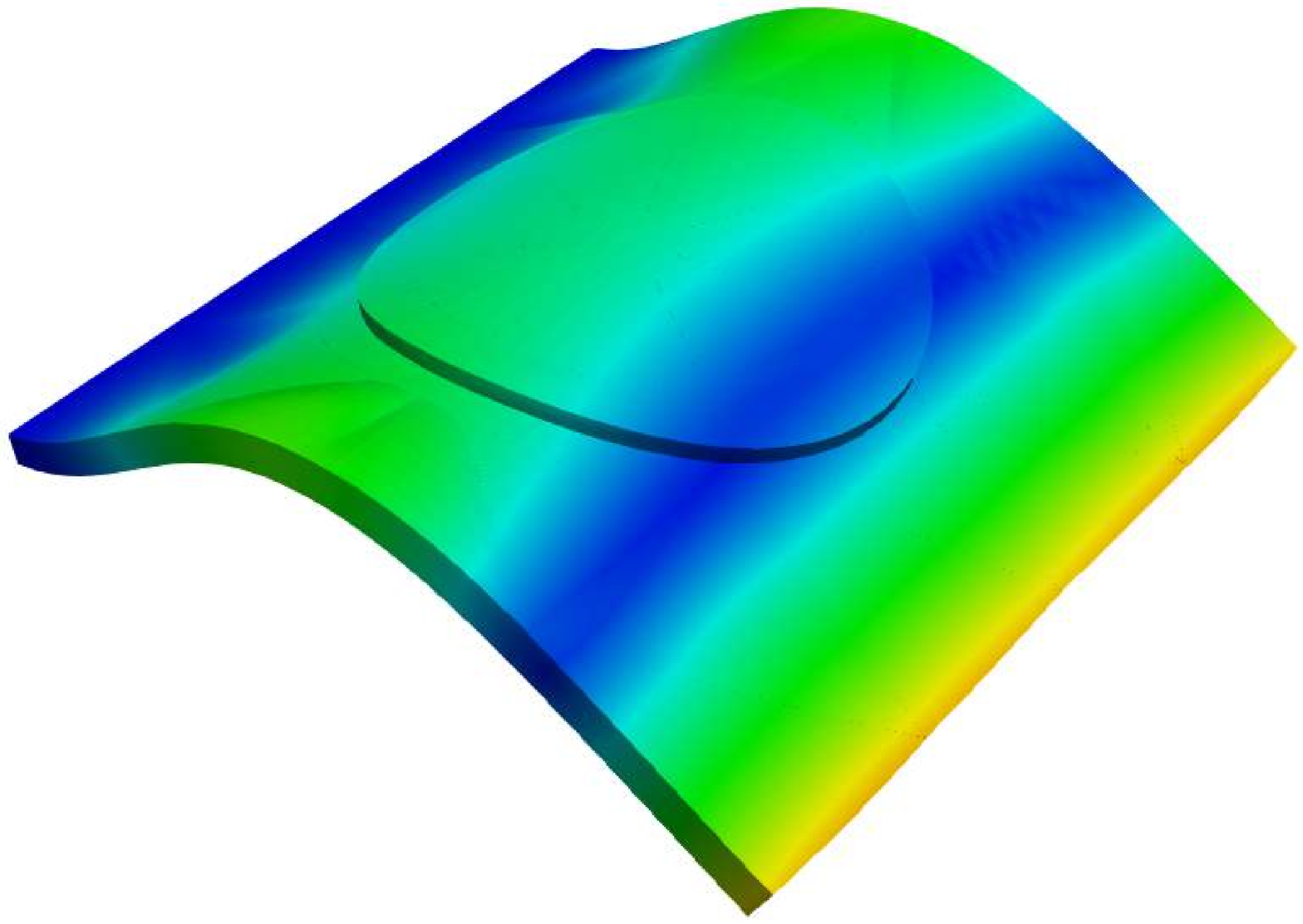}
%		\subcaption{Thrid Eigenform}
%		\label{fig:circPatch_EF3}
	\end{subfigure}

	\caption{Contour plot of the first three eigenforms of plate with circular patch actor}
	\label{fig:circPatch_EF}
\end{figure}

Finally we show results of the eigenvalue calculation. The lowest five eigenfrequencies of the assembly for open ($OC$) and short circuit ($SC$) configuration are listed in Table \ref{tab:circPatch_EF}. 
 For the short cut configuration reference eigenfrequencies were calculated in \textit{ABAQUS}.  The relative frequency error $\Delta f_{SC}=\frac{f_{SC,ref}}{f_{SC}}-1 $  is listed in Table \ref{tab:circPatch_EF}. The lowest three eigenforms are plotted in Figure \ref{fig:circPatch_EF}.

%\begin{comment}
\begin{table}[htp]
	\centering

	\caption{Eigenfrequecies for open circuit $f_{OC}$ and shortcut $f_{SC}$ }
	\label{tab:circPatch_EF}
	\begin{tabular}{|c||c|c||c|c|}
		\hline
		$\#$ & $f_{OC} / \operatorname{kHz} $ & $f_{SC} / \operatorname{kHz}$& $f_{SC,ref} / \operatorname{kHz}$ & $\Delta f_{SC} $  \\ \hline \hline
		1  & 1.2681  & 1.2650  &1.2647  & -0.001066   \\ \hline
		2  & 3.8039  & 3.8039  &3.7987  & -0.000324  \\ \hline
		3  & 8.9729  & 8.8098  &8.8475  & 0.003226  \\ \hline
		4  & 11.7410 & 11.5219 &11.6111 & 0.012628 	 \\ \hline
		5  & 12.0651 & 12.0646 &12.1416 & 0.011402 	 \\ \hline
	\end{tabular}
\end{table}
%\end{comment}

%\begin{comment}
%\end{comment}

\subsection{Radially polarized semicylinder}
\begin{figure}[htp]
	\centering
	\begin{subfigure}{0.4\textwidth}
		
		\psfrag{x}{$x$}
		\psfrag{y}{$y$}
		\psfrag{z}{$z$}
		\psfrag{clamped}{\footnotesize {clamped face}}
		\psfrag{tp}[l]{\footnotesize {tip displacement}}
		\psfrag{phi100}[r]{$\phi=100V$}
		\psfrag{phi0}{$\phi=0V$}
		\psfrag{ra}[r]{\begin{tabular}{@{}r@{}}
				\footnotesize {radially polarized}\\
				\footnotesize {piezoelectric material}
		\end{tabular}}
		\psfrag{dm}[r]{$r=15~mm$}
		\psfrag{hp}{$1~mm$}
		\includegraphics[width=0.75\textwidth]{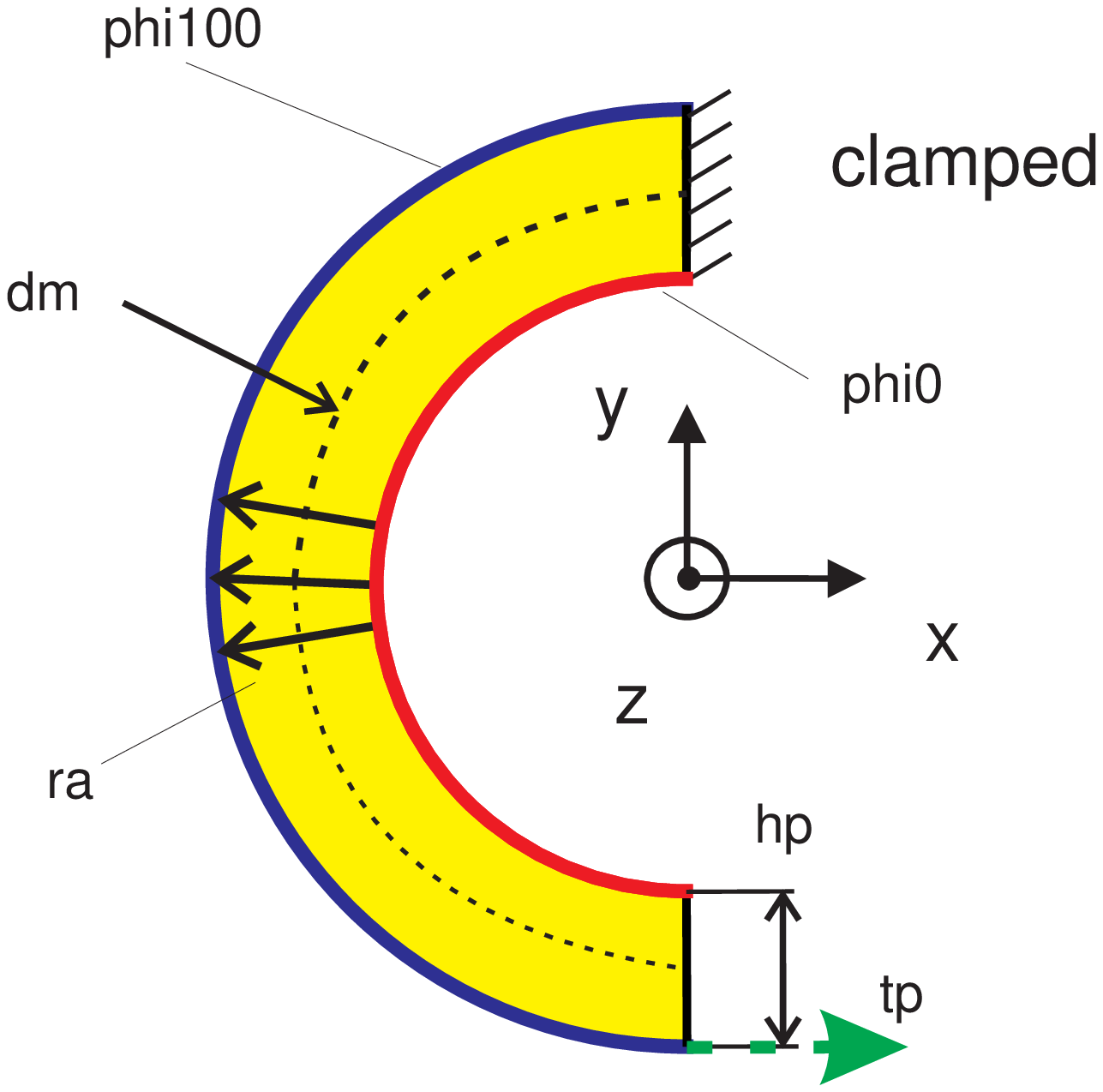}
		\caption{Sketch }
		\label{fig:scetch_semizylinder}
	\end{subfigure}
	\begin{subfigure}{0.3\textwidth}
		
		\includegraphics[width=0.9\textwidth]{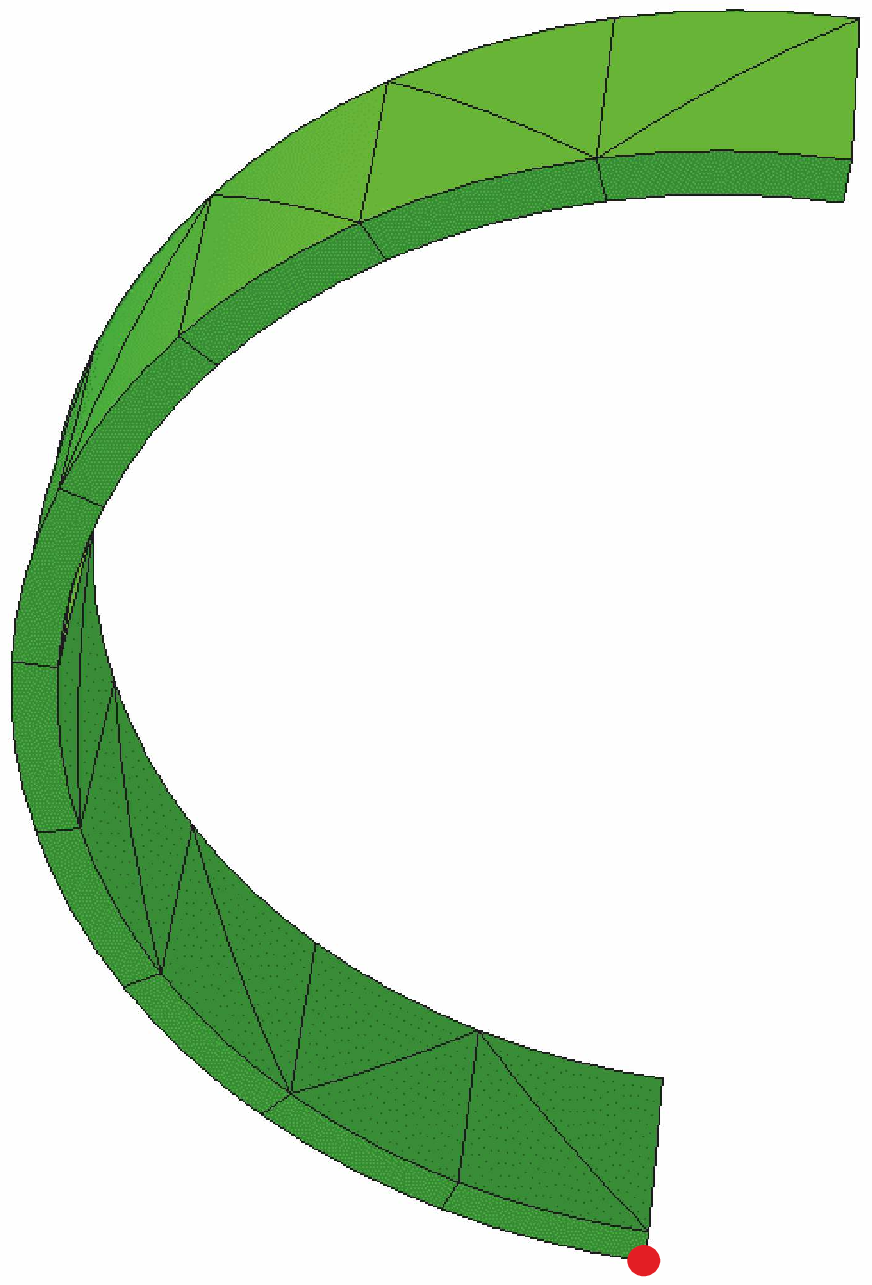}
		\subcaption{Mesh (mesh size $h=5~\operatorname{mm}$)}
		\label{fig:mesh_semicly}
	\end{subfigure}
	\caption{Sketch of radially polarized piezoelectric semicylinder and mesh used in  \textit{Netgen/NGSolve} }
	\label{fig:semicyl}
\end{figure}

Our second example is a radially polarized piezoelectric semicylinder that was first presented in \cite{ZOUARI:2010}. The semicylinder is of   diameter  $15\operatorname{mm}$, thickness  $1\operatorname{mm}$ and length (in $z$ direction)  $5\operatorname{mm}$. A sketch of the geometry is shown in Figure \ref{fig:semicyl}.  We re-use the material parameters from Table \ref{tab:mat_piezp} in the local frame of the polarization direction. To the electrodes, located at the inner and outer cylinder surfaces, we apply a potential difference of $100~V$.  %This example was originally presented by Zouari et.al \cite{ZOUARI:2010}.
%Due to  gradient of the displacement across the thickness caused by the longitudinal $d_{33}$ effect bending stress is caused \cite{ZOUARI:2010}.
\par 

%\begin{comment}
%\end{comment}

We compute a reference solution in \textit{ABAQUS} 6.14 using $10540$ C3D20E-elements (four across the thickness, 17 across the length and 155 across the circumference), in total $210k$ degrees of freedom. The reference tip displacement at the marked corner point in Figure \ref{fig:semicyl}(\subref{fig:mesh_semicly}) is $u^{tip}_{ref}=2.75525\operatorname{\mu m}$.

\begin{comment}
We use the variational formulation (\ref{eq:TDNNS_update})-(\ref{eq:TDNNS_update_rhs}) to implement the constitutive equations. The deformation gradient for this example is given by
\begin{equation}\label{eq:def_grad_semicyl}
	F_M=F_{M,1}\times F_{M,2}= \begin{bmatrix}
			cp & sp & 0\\
			-sp & cp & 0\\
			0 & 0 & 1 
	   \end{bmatrix} \cdot
	   \begin{bmatrix}
	   		0 & 0 & 1 \\ 0& 1 & 0 \\ 1&0 &0
	   \end{bmatrix}.
\end{equation}
with $cp=\frac{x}{\sqrt{x^2+y^2}}$ and $sp=\frac{y}{\sqrt{x^2+y^2}}$.% $cp$ and $sp$ are sine and cosine of angle between material and global coordinate system. 
Therefore $F_{M,1}$ performs a rotation in $z$,  $F_{M,2}$ performs an additions rotation in $y$ by an angle of $\pi/2$, which can as well be interpreted as a change of $x$ and $z$ axis. 
\end{comment}

The mesh used for calculations in \textit{Netgen/NGSolve} is shown in Figure \ref{fig:mesh_semicly}. Again a very coarse prismatic mesh consisting of  one element along the length of semicylinder and 10 elements along the circumference is used. This corresponds to a mesh size $h=5$. A second (finer) mesh with two elements in z-direction and 20 elements along the circumference (mesh size $h=2.5~\operatorname{mm}$) is studied. For both meshes  we use only one element in thickness direction. We performed calculations with different orders of interpolation $k$ for displacements and stresses, interpolation order for the electric potential is $k+1$. 
For each calculation the tip displacement $u^{tip}$ at the corner point is evaluated. We evaluate the relative difference of results defined as $\delta_{rel}=u^{tip}/u^{tip}_{ref}-1 $. The results of this convergence study are shown in Table \ref{tab:semicyl} together with the number of degrees of freedom $n_{dof}$ used in \textit{Netgen/NGSolve}. 
%A small convergence study is shown in Table \ref{tab:semicyl}. 
We see that the accuracy of the solution increases both, when decreasing the mesh size or increasing the polynomial order.
More detailed convergence studies can be found in \cite{PechsteinMeindlhumerHumer:2018}.

%\begin{comment}
\begin{table}[htp]
	\centering
	\caption{Relative difference of tip displacement computed by \textit{Netgen/NGSolve} and \textit{ABAQUS}}
	\begin{tabular}{|c||c| c||c|c|}
		\hline
		  & \multicolumn{2}{|c||}{$h=5$} & \multicolumn{2}{c|}{$h=2.5$} \\ \hline
		  k & $n_{dof}$ & $\delta_{rel}$ & $n_{dof}$ & $\delta_{rel}$ \\ \hline
		  1 & 845 & $-0.016034$ & 2801 & $-0.002189$ \\ \hline
		  2 & 2209 & $0.006113$ & 7543 & $0.000552$ \\ \hline
		  3 & 4521 & $0.001890$ & 15729 & $0.000320$ \\ \hline
	\end{tabular}	
	\label{tab:semicyl}
\end{table}

\section{Conclusion} In our contribution we have provided prismatic and hexahedral elements of arbitrary  polynomial order for the discretization of flat possibly curved piezoelectric geometries by the TDNNS-method. 
The inverse iteration has been proposed to calculate eigenvalues and -forms. %\red{and it was shown why it can be used?}

In our numerical examples we have shown, that accurate simulation of flat curved piezoelectric structures can be done with only one element in thickness direction. % We have seen a significant reduction of degrees of freedom required for simulations.
 Appropriate results were achieved for displacements and stresses as well as for frequencies, while the number of degrees of freedom could be reduced significantly. 

\section*{Acknowledgment}
Martin Meindlhumer acknowledges support of Johannes Kepler University Linz, Linz Institute of Technology (LIT).

%\paragraph{Wordcount:} 4860 words

\bibliographystyle{unsrt}
\bibliography{Smart} 
\clearpage

\end{document}